\documentclass[leqno]{article}

\usepackage[letterpaper,top=2cm,bottom=2cm,left=3cm,right=3cm, marginparwidth=2cm,
]{geometry}
\usepackage[ansinew]{inputenc}
\usepackage{lineno,hyperref}
\modulolinenumbers[5]

\usepackage[usenames, dvipsnames]{color}

\usepackage{amsthm}
\usepackage{graphicx}
\usepackage{microtype}
\usepackage[usenames,dvipsnames]{xcolor}
\usepackage{amsmath,amssymb,amsfonts,latexsym,cancel}
\usepackage{amsmath, amstext, amsfonts}
\usepackage{float}
\usepackage{caption}
\usepackage{subcaption}
\usepackage{appendix}
\usepackage{multirow}

\newcommand{\F}{\mathbb{F}}
\newcommand{\B}{\mathbb{B}}

\usepackage{url}

\newtheorem{theorem}{Theorem}[section]

\newtheorem{Definition}[theorem]{Definition}

\begin{document}

{\centering\LARGE Generation of dynamical S-boxes via lag time chaotic series for cryptosystems\par}
\bigskip

{\centering\large B. B.~Cassal-Quiroga$^1$ and E.~Campos-Cant\'on$^2$ \par \footnotetext[2]{Corresponding Author}}
{\centering\itshape Divisi\'on de Matem\'aticas Aplicadas, Instituto Potosino de Investigaci\'on Cient\'ifica y Tecnol\'ogica A. C., Camino a la Presa San Jos\'e 2055, Col. Lomas 4 Secci\'on, C.P. 78216, San Luis Potos\'i, S.L.P., M\'exico. $^1$bahia.cassal@ipicyt.edu.mx, $^2$eric.campos@ipicyt.edu.mx\par}

\begin{abstract}
In this work, we present an algorithm  for the design of $n\times n$-bits  substitution boxes (S-boxes) based on time series of a discrete dynamical system with chaotic behavior. The elements of a $n\times n$-bits  substitution box are given by binary sequences generated by time series with uniform distribution. Particularly, time series with uniform distribution are generated via two lag time chaotic series of the logistic map. The aim of using these two lag time sequences is to hide the map used, and thus U-shape distribution of the logistic map is avoided and uncorrelated S-box elements are obtained. 
The algorithm proposed is simple and guarantees the generation of S-boxes, which are the main component in block cipher, fulfill the strong S-box criteria: bijective; nonlinearity; strict avalanche criterion; output bits independence criterion; criterion of equiprobable input/output XOR distribution and; maximum expected linear probability. 
The S-boxes that fulfill these criteria are commonly known as ``good S-boxes". 
Finally, an application based on polyalphabetic ciphers principle is developed to obtain uniform distribution of the plaintext via dynamical S-boxes.\\

\bf{S-box; block cipher; dynamical S-Box; chaos; lag time chaotic series.}
\end{abstract}

\section{\label{sec:Introduction}Introduction}
\label{introduction}

Nowadays, we are in the era of informatics and due to the large number of attacks, it is important to adequately protect the information to avoid possible misuse of it. The aforementioned comment motivates the generation of different approaches to have secure cryptographic systems. In general cryptosystems can be divided in two classes: stream cipher and block cipher. The stream cipher combines bit to bit, the sequences of bits generated by plaintext and pseudorandom numbers. The block cipher takes blocks of plaintext which are encrypted by substitution using an S-box and cyclic shifting. The substitution box (S-box) is the core component of block cipher. The S-boxes give the cryptosystems the confusion property described by Shannon \cite{Shannon1949}, used in conventional block ciphers such as the Data Encryption Standard (DES) and the Advanced Encryption Standard (AES). In these cryptosystems the security depends mainly on the S-box properties that are used. The criteria that a strong S-box fulfill, also known as ``good S-boxes",  are: bijection; nonlinearity; strict avalanche criterion (SAC); the output bit independence criterion (BIC) \cite{Adams1990structured}. Other desirable characteristics are to be resistant to linear and differential  cryptanalysis attacks. The construction of cryptographically secure S-boxes is a field of interest in the cryptography area.

In recent years, many papers have been reported and are focused on studying cryptosystems based on chaos \cite{jakimoski2001chaos, CHEN, WANG, LAMBIC, belazi2017efficient, OZKAYNAK, Liu, ccavucsouglu2017novel, Guesmi, Tian}, this is, because of the relationship that exists between the chaotic system properties and the cryptosystem properties. In \cite{ALVAREZ} the relationship between these properties are given, for instance, confusion is related with ergodicity, the diffusion property with sensitivity to initial conditions and the deterministic dynamic with the deterministic pseudo-randomness. Taking advantage of the properties of chaotic systems, we propose a  strong and dynamic S-Box that complies with the different criteria of good S-boxes.

Regarding the generation of S-box based on chaos, some algorithms have been developed
using discrete dynamical systems. For example, in \cite{jakimoski2001chaos, CHEN, WANG, LAMBIC, belazi2017efficient}, the generation of substitution boxes were introduced through a single time series of a map or by combining two time series of different maps. Nevertheless, these algorithms do not guarantee that the series used have a uniform distribution, as in our approach based on two lag time chaotic series derived from the logistic map. In the same way, there are algorithms based on continuous chaotic dynamical systems \cite{OZKAYNAK,Liu, ccavucsouglu2017novel}. Also there are algorithms based on the mixing of time series of continuous and discrete dynamical systems \cite{Guesmi,Tian} and in \cite{Azkaynak2013} the algorithm is build via time-delay series. 
The advantage of using discrete chaotic dynamical systems is that from one iteration to another the elements of the time series are de-correlated, however this does not happen if a continuous chaotic dynamical system is used, the elements of the time series are strong  correlated. Therefore, many iterations are needed and the calculation of the mutual information between elements of the time series is necessary to be able to say when they are de-correlated, which implies higher computational cost. 

In chaos-based encryption schemes, pseudo-random sequences based on chao\-tic maps are generally used as one time pad for encrypting messages. Since encryption schemes, based on low dimensional chaotic map, have low computational complexity, they can be analyzed with low computational cost using iteration and correlation functions \cite{ZHOU}.
Time-delay chaotic series have complex behavior and erase the trace of the mapping that generates them. Using these aforementioned time series, S-boxes can be designed and provide better nonlinearity criterion, that ensure good statistical properties in the generators. 

In this paper, a method to obtain dynamical good S-boxes is presented based on the generation of lag time series from the logistic map. Using this kind of lag time series, it is possible to generate S-boxes that have the capability of hiding the map used to build them. With this approach, based on these lag time series, pseudo-random series are generated with good statistical properties, more details can be found in \cite{Garcia2014Campos}. This novel algorithm for S-box generation is based on cryptographically secure pseudo-random number generator. The rest of the paper is organized as follows: In Section \ref{Criterias}, the criteria for a ``good'' $n \times n$ bit S-box are described. In Section \ref{Sec secuencias}, it is presented a dynamical analysis of logistics map. In Section \ref{Sec algorithm}, the proposed scheme to generate dynamical S-box based on pseudo-random bit generator is presented. In  Section \ref{Performance test of S-box}, the performance analysis of an obtained S-box and its comparison with other S-box reported in the literature is given. An application of the obtained S-boxes to hide an image is presented in Section \ref{Application of dynamical S-boxes}. Finally conclusions are drawn in Section \ref{conclusions}.

\section{Criteria for a good n$\times$n bit S-box} 
\label{Criterias}

A collection of six criteria reported in the literature for generate cryptographically good S-boxes has been made. These criteria  are: bijective; nonlinearity; strict avalanche criterion; output bits independence criterion, equiprobable input/output XOR distribution; and maximum expected linear probability. Before addressing these properties it is necessary to give some preliminaries about Boolean functions.

Let $\B=\{0,1\}$ be a binary set which is endowed with two binary operations, called addition (denoted by $\oplus$ XOR operation) and multiplication (denoted by~$\cdot$ AND operation). Let $(\B,\oplus,\cdot)$ be a field which will be denoted by $\F$, where the binary operations are given by the table \ref{tabla1}.
\begin{table}[h!]
	\centering
	\begin{tabular}{|l|l|l|}
		\hline
		$\oplus$ &0&1\\
		\hline
		0&0&1\\
		\hline
		1&1&0\\
		\hline
	\end{tabular}
	\; \; \;  
	\begin{tabular}{|l|l|l|}
		\hline
		$\cdot$ &0&1\\
		\hline
		0&0&0\\
		\hline
		1&0&1\\
		\hline
	\end{tabular}
	\caption{Addition and multiplication binary operations.}
	\label{tabla1}
\end{table}

An $n\times n$ S-box is a vectorial Boolean function $S : \F^n \to \F^n$, where $\F^n$ is a vectorial space, and  is defined as:
\begin{equation}
S(x) = (f_1(x),f_2(x),\cdots , f_n(x)),
\end{equation} 
where $x = (x_1, x_2, \cdots , x_n)^\top \in \F^n$ and each of $f_i'$s for $1 \leq i \leq n$ is a Boolean function.
A Boolean function is a mapping $f : \F^n \to \F $ by considering all inputs in $f$, $f_i$ can be seen as a column vector of $2^n$ elements. The functions $f_i'$s are component functions of $S$.

Some basic definitions can be found in \cite{Dimitris}. 
\begin{Definition}
	A Boolean function with algebraic expression, where the degree is at most one is called an {\bf affine Boolean function}. The general form for n-variable affine function is:
	$$f_{affine}(x_1, x_2, x_3, \ldots , x_n) = w_n\cdot  x_n \oplus w_{n-1}\cdot x_{n-1} \oplus \ldots  \oplus w_2 \cdot x_2 \oplus w_1 \cdot x_1 \oplus w_0,$$
	where $w_i\in \B$ are coefficients, and $x_i\in \B$ are variables, with $i=0,1,\ldots,n$.
\end{Definition}

\begin{Definition}
	A {\bf linear Boolean function} is defined as follows
	$$L_w(x)=w_n\cdot x_n \oplus w_{n-1}\cdot x_{n-1} \oplus \ldots \oplus w_1\cdot x_1,$$ 
	where $x_i,w_i\in\B$, with $i=1,\ldots,n$.
\end{Definition}

The set of affine Boolean functions is comprised by the set of linear Boolean functions and their complements,{\it i.e.}, all functions of the form
$$ A_{w,c}(x)= L_w(x)\oplus c,\; c\in \B.$$

A useful representation of a Boolean function $f_i$, with $i=1,\ldots,n$, is given by the polarity truth table defined as follows.

\begin{Definition}
	A {\bf polarity truth table} is defined as follows
	$$\hat{f} (x)=(-1)^{f(x)},$$
	where $x\in \F^n$,  $\hat{f}$ maps the output values of the Boolean function from the set $\{0, 1\}$ to the set $\{ -1, 1\}$, $i.e.$,
	$$\hat{f}: \{0, 1\} \to \{ -1, 1\}.$$
	A linear Boolean function in polarity form is denoted as $\hat{L}_w(x)$.
\end{Definition}

\begin{Definition}
	The {\bf Walsh Hadamard transform (WHT)} of a Boolean function $f$ is defined as
	$$\hat{F}_f(w)=\sum_{x\in \B^n}{\hat{f}(x)\hat{L}_w(x)}.$$
\end{Definition}

The WHT measures the correlation between the Boolean function $f$ and the linear Boolean function $\hat{L}_w$ with $x \in \F^n$. 

\subsection{Bijective Criterion }
\label{Bijective}

Let $S(x)$ be an S-box, which is bijective if and only if their Boolean functions $f_i$ satisfy the following condition:  
\begin{equation}\label{bijective}
wt(a_1 \cdot f_1 \oplus a_2\cdot f_2 \oplus \cdots \oplus a_n\cdot f_n)=2^{n-1},
\end{equation}

where $a_i \in \F$, $(a_1,a_2, \cdots ,a_n)\ne(0, 0, \cdots , 0)$ and $wt(\cdot)$ is the Hamming weight \cite{Adams1990structured, adams1989good}, the corresponding S-box is guaranteed to be bijective.

\subsection{Nonlinearity criterion}
\label{Nonlinearity}

\begin{Definition} \cite{Tian2017Chaos}
	The nonlinearity of a Boolean function $f : \F^n \to \F $ is denoted
	by 
	\begin{equation}
	\label{distance Nonlinearity}
	N_f = \min_{l\in A_{w,c}(x)} d_H( f , l),
	\end{equation}
	where $A_{w,c}(x)$ is an affine function set, $d_H( f , l)$ is the Hamming distance between $f$ and $l$.
\end{Definition}

The minimum distance between two Boolean functions can be described by means of the Walsh spectrum \cite{Millan1998}: 
\begin{equation}\label{nonlinearity equation}
\min_{l\in A_{w,c}(x)} d_H( f , l)  =2^{n-1} (1-2^{-n} \max_{\omega \in \F^n} |\hat{S}_{(f)}(\omega)|),
\end{equation}

where the Walsh spectrum of $f(x)$ is defined as follows:
\begin{equation}\label{spectrum nonlinearity}
\hat{S}_{(f)}(\omega)=  |\hat{F}_f(w)|_{w\in\B^n} =\sum_{x\in \F^n}(-1)^{f(x)\oplus x \bullet \omega},
\end{equation}
with $\omega \in \F^n$ and $x\bullet \omega$ is the dot product between $x$ and $
\omega$ as:
\begin{equation}\label{dot product nonlinearity}
x\bullet \omega = x_1 \cdot\omega_1 \oplus \cdots \oplus x_n \cdot\omega_n.
\end{equation}

\subsection{Strict Avalanche Criterion (SAC)} 
\label{SAC}

This criterion was first introduced by Webster and Tavares \cite{Webster1986}. A Boolean function $f$ satisfies SAC if complementing any single input bit changes the output bit with the probability one half. So, more formally, a Boolean function $f$ satisfies SAC, if and only if
\begin{equation}\label{SAC equation}
\sum_{x\in \F^n}f(x)\oplus f(x\oplus e_i)=2^{n-1},    ~~~\forall i:1 \leq i \leq n,
\end{equation}
where $e_i\in \F^n$ such that $wt(e_i)=1$.

\subsection{Output Bits Independence Criterion (BIC)}
\label{BIC}
Output Bit Independence Criterion is another desirable criterion for an S-box that should be satisfied, introduced by Webster and Tavares \cite{Webster1986}.
It means that all the avalanche variables should be pairwise independent for a given set of avalanche vectors generated by the complementing of a single plaintext bit.

Adam and Tavares introduced another method to measure the BIC that for the Boolean functions, $f_i$ and $f_j$ $(i\ne j )$ of two output bits in a S-box, if  $f_i\oplus f_j$ is highly nonlinear and come as close as possible to satisfy SAC \cite{Adams1990structured}. Additionally, $f_i\oplus f_j$ can be tested with a Dynamic Distance (DD). The DD of a function $f$ can be defined as:
\begin{equation} \label{dynamic distance equation}
DD(f)=\max_{\substack{d\in \F^n\\ wt(d)=1}}\frac{1}{2}\left| 2^{n-1}-\sum_{x=0}^{2^n -1} f(x)\oplus f(x\oplus d)\right|.
\end{equation}

If the value of DD is a small integer and close to zero, the function $f$ satisfies the SAC. 

\subsection{Criterion of equiprobable Input/Output XOR Distribution}
\label{MEDP} 

Biham and Shamir \cite{biham1991differential} introduced differential cryptanalysis which attacks S-boxes faster than brute-force attack. It is desirable for an S-box to have differential uniformity. This can be measured by the maximum expected differential probability (MEDP). Differential probability for a given map $S$ can be calculated by measuring differential resistance and is defined as follows:
\begin{equation} \label{MEDP equation}
DP_f=\max_{\Delta x\ne 0, \Delta y} \left( \frac{ \# \{ x \in \F^n| S(x) \oplus S(x \oplus \Delta x)= \Delta y \} }{2^m}\right),
\end{equation}

\noindent where $2^n$ is the cardinality of all the possible input values ($x$), $\Delta x$ and $\Delta y$ are called input and output differences, respectively, for the $S$. Thus, the smaller value of $DP_f$ gives better cryptographic property, i.e., its resistance to differential cryptanalysis.

\subsection{Maximum Expected Linear Probability}
\label{MELP}

The Maximum Expected Linear Probability (MELP) is the maximum value of the unbalance of an event. Given two randomly selected masks $a$ and $b$, and $a$ is used to calculate the mask of all possible values of an input $x$, and use $b$ to calculate the mask of the output values of the corresponding S-box. The parity of the input bits mask $a$ is equal to the parity of the output bits the mask $b$. MELP of a given S-box can be computed by the following equation:
\begin{equation}\label{MELP equation}
LP_f=\max_{a,b \in \F^n \backslash\{0\} } \left( 2^{-n} \sum_{x\in \F^n}(-1)^{a \cdot x + b \cdot f(x)}\right)^2.
\end{equation}

The closer is MELP to zero, the higher the resistance against linear cryptanalysis attack will be.

\section{Analysis of the Logistic map}
\label{Sec secuencias}

The logistic map is a discrete-time demographic model analogous to the logistic equation first created by Pierre François Verhulst, which is described by the following differential equation
$$ \frac{dx}{dt}=rx\left(1-\frac{x}{K}\right), $$
where $x$ is the state variable of the system, $r$ is a parameter related with the rate of maximum population growth and $K$ is the so-called carrying capacity (i.e., the maximum sustainable population). So $x\leq K$, when $x=K$ the population stops growing. 
Robert May \cite{may1976simple} popularized this differential equation to one of the most famous discrete dynamical systems, the logistic map, which is defined as follows:
\begin{equation} \label{logistic map}
f_\alpha (x_i)=\alpha x_i (1-x_i), 
\end{equation}

\noindent where $x$ is the state variable of the logistic map and, $\alpha$ is the only  parameter of the system instead of two as its analogous continuous model. The use of a single parameter was possible because the logistic map was normalized, {\it i.e.}, $f_\alpha:[0,1]\to [0,1]$, for the bifurcation parameter $ \alpha \in \lbrack 0,4\rbrack$ and $x_0\in[0,1]$. Nevertheless, in the context of mathematics, the values of the parameter $ \alpha$ are not restricted to the interval $\lbrack 0,4\rbrack$, so mathematically, it is possible to consider negative values \cite{dendrinos1993socio}. As mentioned above, the logistic map is now studied in the interval $\lbrack -2,0 )$ for cryptographic purposes. Now, we study the mapping behavior in the two intervals and they assured with $\alpha \in \lbrack-2,4\rbrack$ the orbits do not scape to infinity for some initial conditions. 
\begin{figure}[h]
	\centering
	\includegraphics[width=0.6\linewidth]{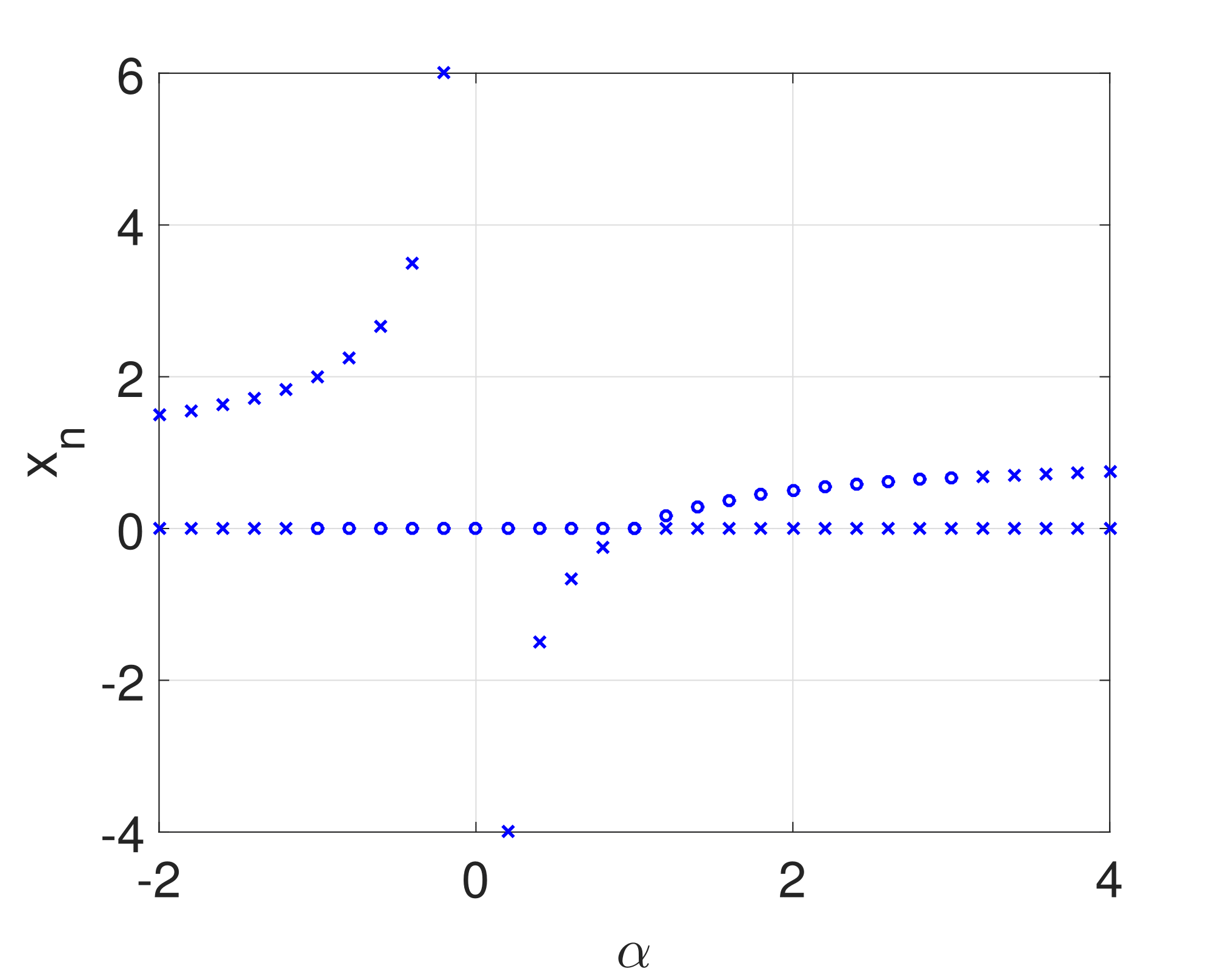}
	\caption{Fixed points stability where an asterisk and a circle denote repulsive and attracting fixed points, respectively.}
	\label{fig:estabilidaddelosptosfijos}
\end{figure}

The dynamical system \eqref{logistic map} presents one or two fixed points located at $x^*_1= 0 $ and at $x^*_2=\frac{\alpha-1}{\alpha}$, for $\alpha \neq 0$. Figure \ref{fig:estabilidaddelosptosfijos}~ depicts the stability of the fixed points where an asterisk and a circle denote  repulsive and attracting fixed points, respectively. These fixed points change their stability according to the parameter $\alpha$, {\it i.e.}, when $|f'(x^*_1)|<1$ and $|f'(x^*_2)|<1$ then the fixed points $x^*_1$ and $x^*_2$ are stable, respectively, and they are unstable when $|f'(x^*_1)|>1$ and $|f'(x^*_2)|>1$. We are interested in the last case because the system presents complex behavior, this is, both fixed points are repulsive, $|f'(x^*_1)|=|\alpha|>1$ and $|f'(x^*_2)|=|-\alpha+2|>1$. The $x^*_1$ fixed point is repulsive for $\alpha<-1$ and $\alpha>1$. On the other hand, 
the $x^*_2$ fixed point is repulsive for $\alpha<1$ but $\alpha\neq 0$, and $\alpha>3$. So the interested values are  $\alpha\in [-2,-1]\cup [3,4]$, this is the condition to have  both repulsive fixed points.

\begin{figure}[h]
	\centering
	\includegraphics[width=0.63\linewidth]{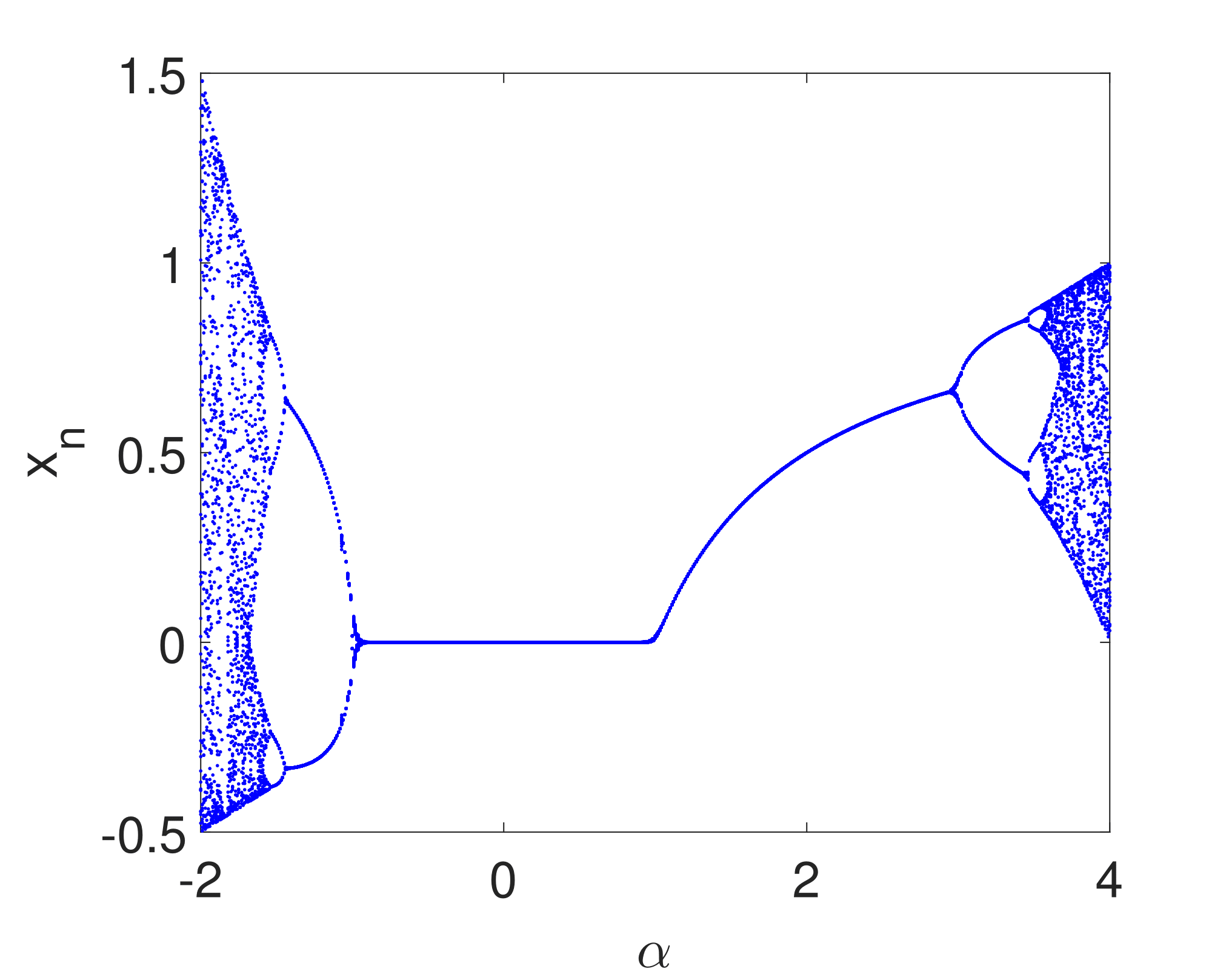}
	\caption{Bifurcation diagram for the logistic map given by Eq. \eqref{logistic map}.  }
	\label{fig:diagramabifurcacion-24}
\end{figure}

The dynamical system \eqref{logistic map}  bifurcates when $|f'(x^*_1)|=1$ and $|f'(x^*_2)|=1$, this happens for $x^*_1$ when $\alpha=-1$ or $1$, and for $x^*_2$ the bifurcations values are given by  $\alpha=1$ and $3$. It is possible to analyze the behavior of the system by means of a bifurcation diagram, which is shown in Figure \ref{fig:diagramabifurcacion-24}. This diagram shows orbits as a function of $\alpha$ parameter and the route to chaos are period-doubling bifurcations at $\alpha=3$ and period-halving bifurcations at $\alpha=-1$. There are intervals for the parameter $\alpha$ near to $-2$ and $4$ where the logistic map $f_\alpha(x)$ behaves chaotically. 

\begin{figure}[h!]
	\centering
	\includegraphics[width=0.6\linewidth]{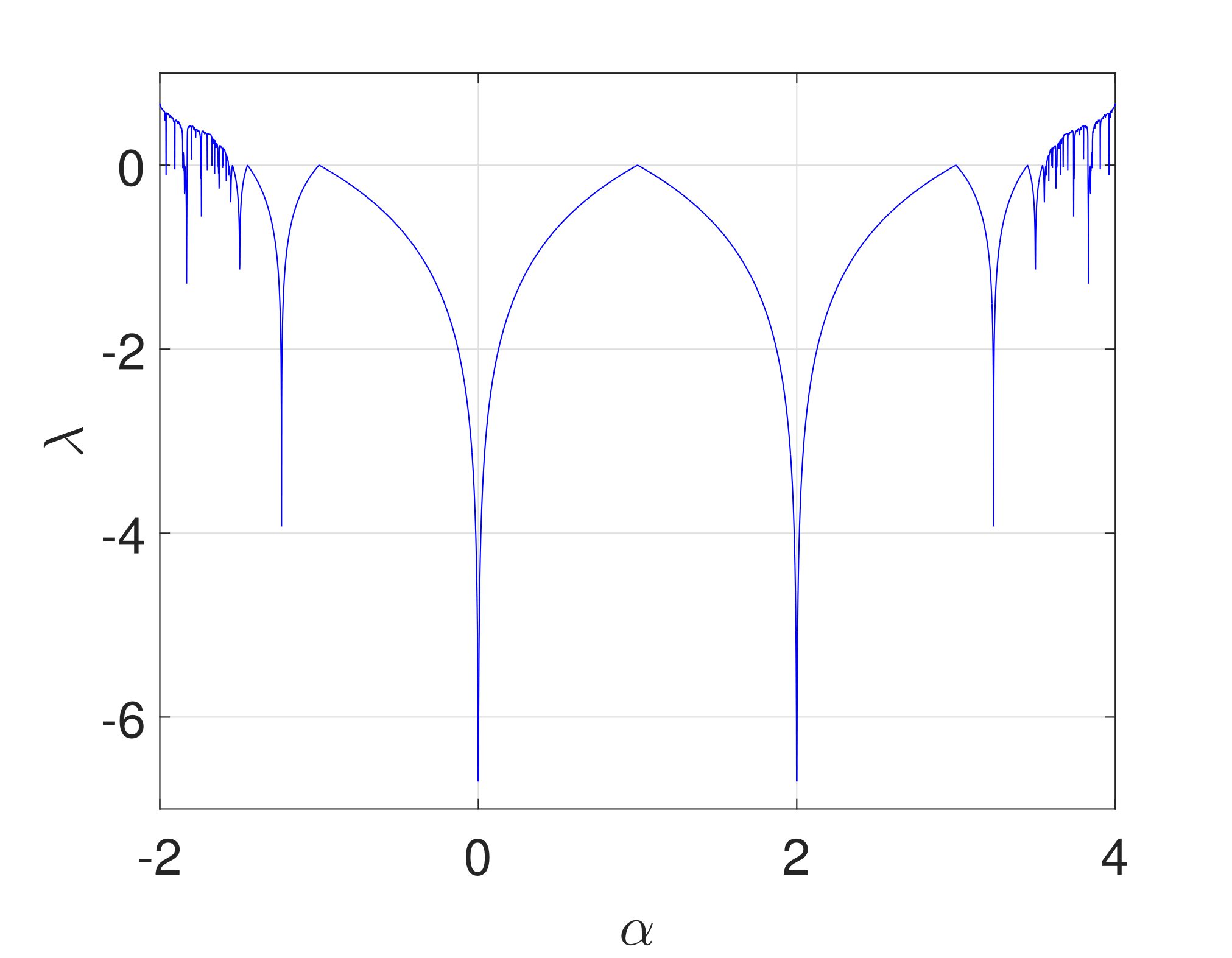}
	\caption{Lyapunov exponent as a function of parameter $\alpha$.}
	\label{fig:exponentedelyapunovmapeologistico1}
\end{figure}

There are several approaches to demonstrate that a system is chaotic, one of them is prove that the dynamical systems fulfills the definition given by Devaney \cite{devaneyintroduction}, other approach is based on the Lyapunov exponent \cite{Lyapunovexponents2004},\cite{Yang2012}. In the same sense, the Lyapunov exponent of Eq. \ref{logistic map} it is shown in Fig. \ref{fig:exponentedelyapunovmapeologistico1}. The graph of Lyapunov exponents is symmetric with respect to $\alpha=1$, the chaotic behavior of the logistic map appears for values of the parameter $\alpha$ near $-2$ and $4$. The local stability of the fixed points are in accordance with the Lyapunov exponent values, for example, when $\alpha\in (-1,3)$ the orbits of the system converge at a fixed point and, when the bifurcations occurs the orbits converge at periodic orbits up to chaos appears. 

\begin{figure}
	\centering
	\includegraphics[width=0.6\linewidth]{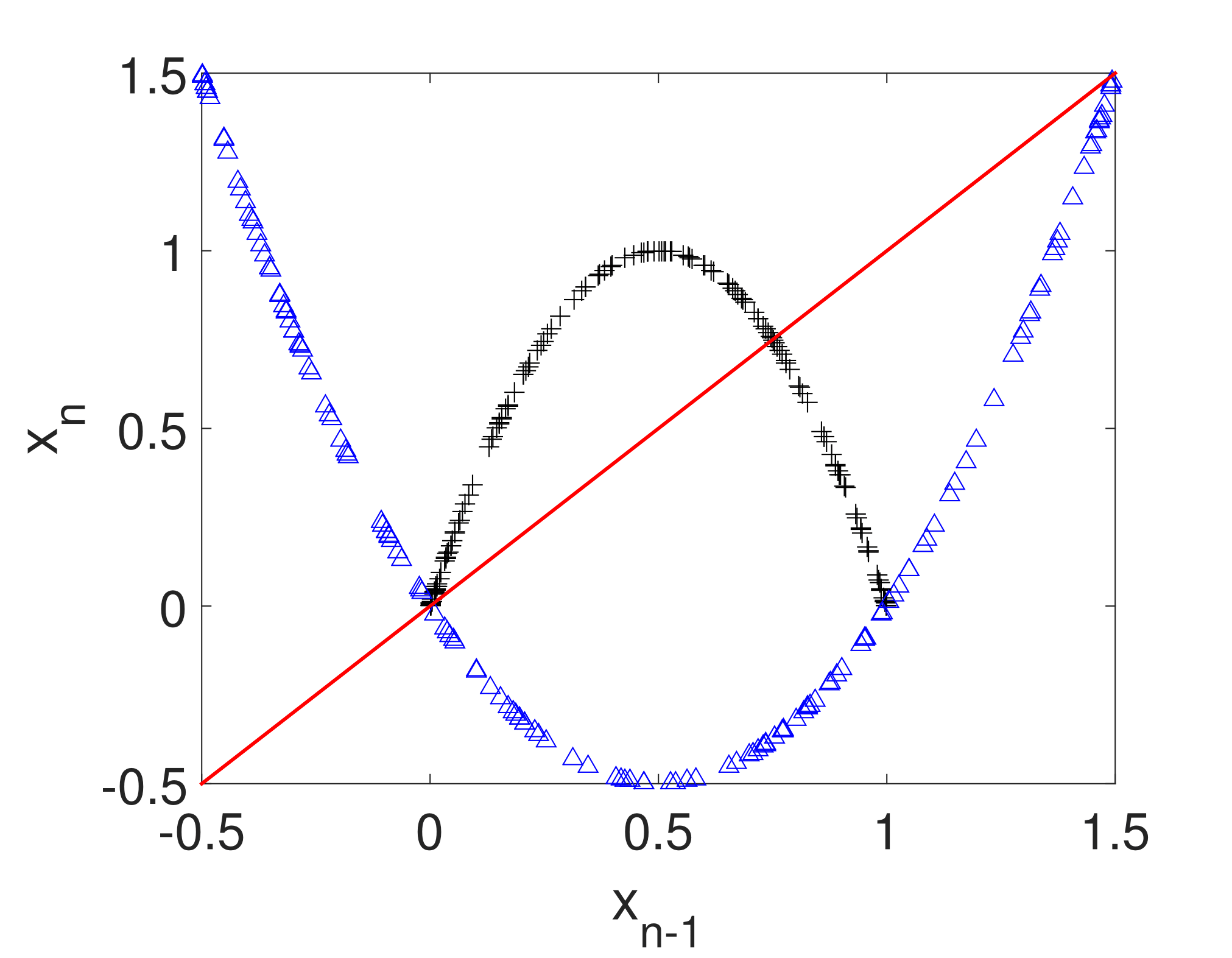}
	\caption{Logistic map for $\alpha=-2$ in blue triangles and for $\alpha=4$ in black crosses.}
	\label{fig:maplogisticoalpa-2y4}
\end{figure}

The aim is to use the logistic map to generate a time series with uniform distribution and without evidencing the mapping used. To achieve this, it is proposed an approach based on two chaotic time series of the logistic map. Based on Lyapunov exponents analysis, the $\alpha$ values are arbitrarily selected within the chaos region, so it is consider $\alpha=-2$ and $4$. In Figure \ref{fig:maplogisticoalpa-2y4} is presented the shape of the logistic map for both parameter values $\alpha=-2$ and $\alpha=4$ in blue triangles and in black crosses, respectively. The logistic map for these parameter values is invariant in different intervals as follows
\begin{equation}
\begin{array}{rl}
f_{-2}&:[-0.5,1.5] \to [-0.5,1.5];\\
f_{4}&:[0,1] \to [0,1].
\end{array}
\end{equation}
It is worth saying that the time series generated with both parameter values have a U-shape distribution. 

\section{Proposed algorithm to generate dynamical S-boxes}
\label{Sec algorithm}

The main idea of the proposed algorithm to generate dynamical S-boxes is based on a Cryptographically Secure Pseudo-Random Number Generator (CSPRNG) via a discrete dynamical system $f_\alpha:I\to I$. Garc\'ia Mart\'inez and Campos Cant\'on \cite{Garcia2014Campos} proposed a CSPRNG using two lag time series generated with the logistic map \eqref{logistic map}.

An orbit $x_0,x_1,x_2,\ldots$ of the logistic map \eqref{logistic map} is defined by giving an initial condition $x_0\in I$.  The interval $I$ is determined by the parameter $\alpha\in \{-2,4\}$, $f_\alpha:I\to I$. Let $M1$ and $M2$ two time series generated with the logistic map by means of the following considerations: {\it i)} given two arbitrary initial conditions $x_{01}$, $x_{02}$, such that, $x_{01}\neq x_{02}$; {\it ii)} two different bifurcation parameters $\alpha_1$ and $\alpha_2$; and {\it iii)} $l$-units of memory for each time series $x_{(i-k_{l-1})1}$, $\ldots$, $x_{(i-k_2)1}, x_{(i-k_1)1}, x_{i1}$ and $x_{(i-k_{l-1})1}$, $\ldots$, $x_{(i-k'_2)2}, x_{(i-k'_1)2}, x_{i2}$. So the orbits have uniform distribution independently the U-shape distribution of the logistic map.
In order to illustrate the algorithm, we have chosen  the  bifurcation parameter values as $\alpha_1=-2$ and $\alpha_2=4$ for the time series $M1$ and $M2$, respectively. These parameter values ensure that the system \eqref{logistic map} has chaotic behavior in both cases. 

To guarantee that the generator will present good statistical properties is necessary to generate time series with uniform distribution and also it is desirable to eliminate the logistic map shape in these new time series. This is achieved by means of the number of lags involved. There are a lot of combinations of delays that are able to de-correlate the shape of the logistic map and the time series, but each delay unit needs memory and processing time. 

\begin{figure}[h!]
	\centering
	\includegraphics[width=0.6\linewidth]{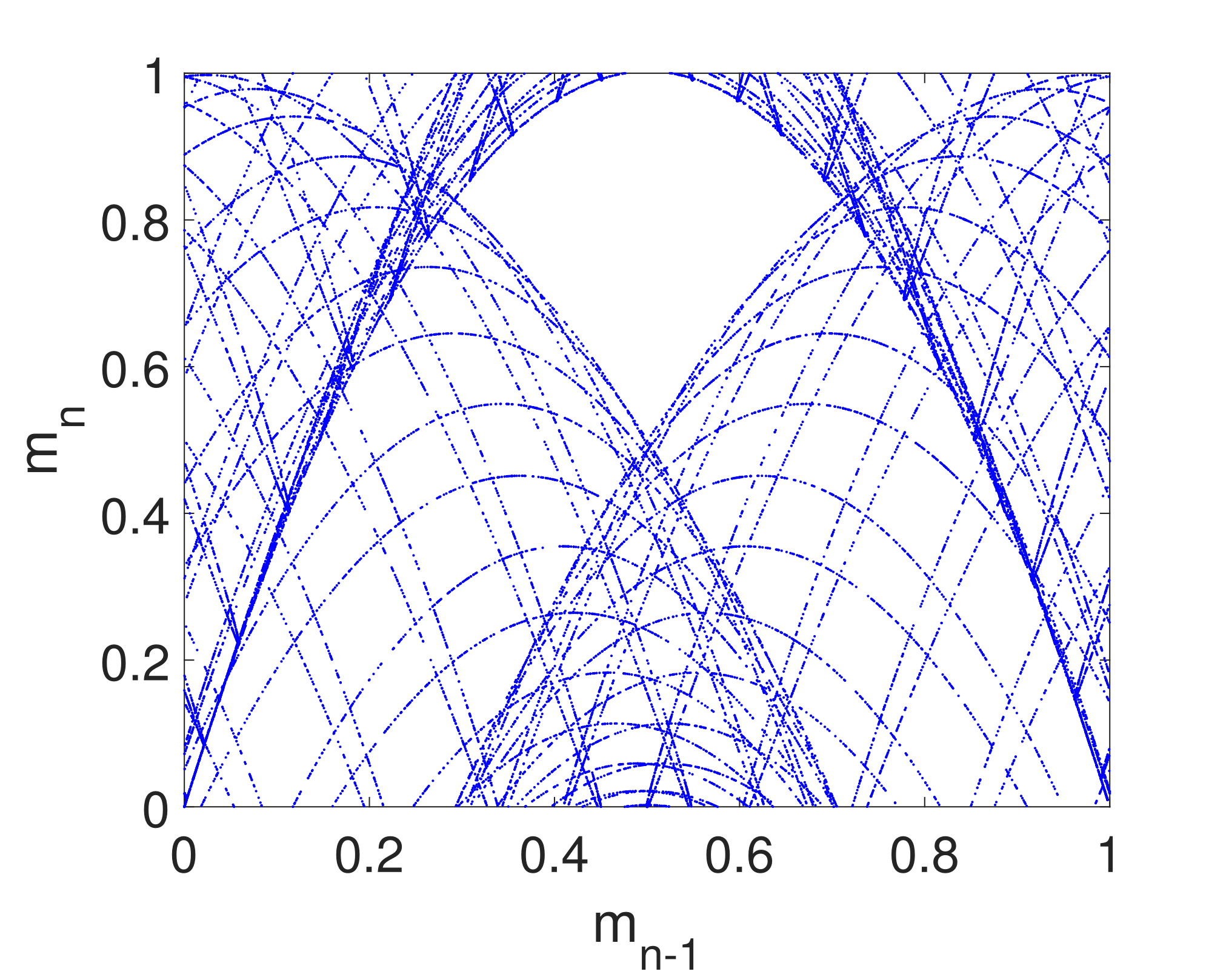}
	\caption{ $(m_{(n-1)2}, m_{n2})$ from the time series $M2(x_{(i-k1)2},x_{i2})$  considering two memory units.}
	\label{fig:retratofasexnconxn12retardos}
\end{figure}

If we analyze time series $M2=m_{02}, m_{12}, m_{22},\ldots$ with two memory units, for $\alpha=4$. The elements $m_{i2}$ of the time series $M2(x_{(i-k1)2},x_{i2})$ are obtained in the following way:
\begin{equation}
\label{2MUTS} 
m_{i2}=M2(x_{(i-k1)2},x_{i2})=x_{(i-k1)2}+x_{i2}, \mbox{mod 1},
\end{equation}
where $k1=5$. In the plot of $ m_{(n-1)2}$ against $m_{n2}$, it is possible to distinguish that the time series $M2$ are generated with the logistic map.
Figure \ref{fig:retratofasexnconxn12retardos} shows $(m_{(n-1)2},m_{n2})$ using two memory units. Because the length of the delay does not matter, the shape of the logistic map always remains, so it is necessary to consider more memory units.

\begin{figure}[h!]
	\centering
	\includegraphics[width=0.6\linewidth]{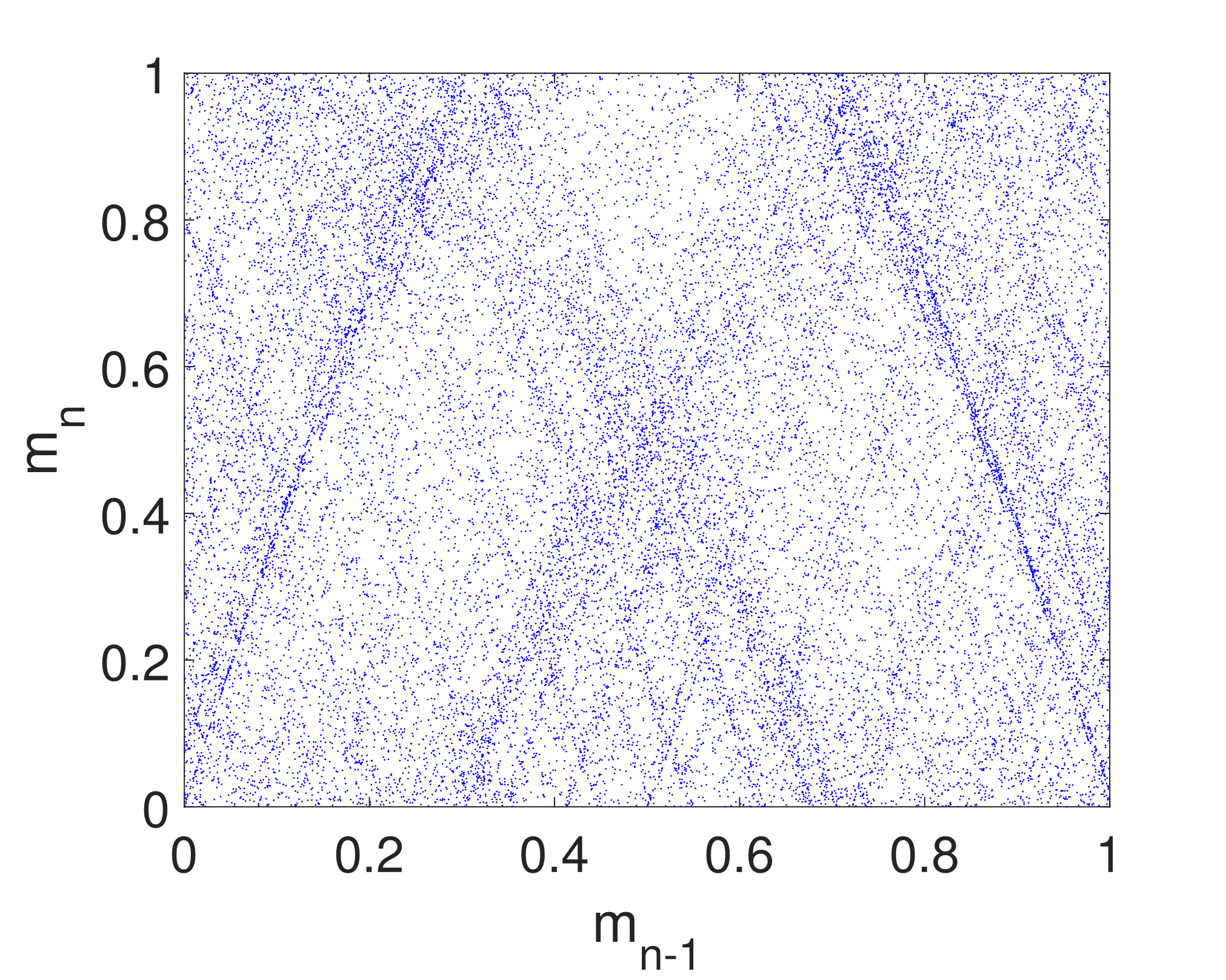}
	\caption{$(m_{(n-1)2}, m_{n2})$ from the time series $M2(x_{(i-k2)2},x_{(i-k1)2},x_{i2})$ considering three memory units.}
	\label{fig:retratofase3retardos}
\end{figure}

Now, if we consider three memory units to obtain the elements of the time series $M2(x_{(i-k2)2},x_{(i-k1)2},x_{i2})$ given as follows
\begin{equation}
\label{3MUTS} 
m_{i2}=M2(x_{(i-k2)2}, x_{(i-k1)2},x_{i2})=x_{(i-k2)2} + x_{(i-k1)2} + x_{i2}, \mbox{mod 1},
\end{equation}
where $k1=10$ and $k2=5$. Now, for this case of three memory units, which are the minimum amount to obtain cloud of points in $(m_{(n-1)2}, m_{n2})$, see Figure~\ref{fig:retratofase3retardos}. The shape of the logistic map almost disappears, so three memory units are enough. The problem with considering more memory units has a computational price of information storage. For this reason two delays $k1$ and $k2$ and the present state of the time series of the logistic map are used. Also the lags must not be contiguous in order to avoid regular patterns which directly affect the test results.

It is considered different delays, {\it i.e.}, $k_2=k'_2=10$, $k_1=6$ and $k'_1=5$ for both time series $M1$ and $M2$. Thus, these series are conformed by the sum of two delay states $x_{(i-10)1}$ and  $x_{(i-5)1}$ with the actual state $x_{i1}$ of the orbit $x_{01},x_{11},x_{21},\ldots$, for $M1$. In the same way for $M2$, $x_{(i-10)2}$,  $x_{(i-6)2}$ and $x_{i2}$ of the orbit $x_{02},x_{12},x_{22},\ldots$.  The values of the time series are limited by the operation mod 1, this guarantees that $M1,M2\in[0,1)\subset \Re$. Explicitly $M1(x_{(i-10)1}, x_{(i-5)1},x_{i1})$ and $M2(x_{(i-10)2}, x_{(i-6)2},x_{i2})$ are expressed in the following way:

\begin{gather}
\label{eqn:M1} 
m_{i1}=M1(x_{(i-10)1}, x_{(i-5)1},x_{i1})=x_{(i-10)1}+x_{(i-5)1}+x_{i1}, \mbox{mod 1},\\
\label{eqn:M2} 
m_{i2}=M2(x_{(i-10)2}, x_{(i-6)2},x_{i2})=x_{(i-10)2}+x_{(i-6)2}+x_{i2}, \mbox{mod 1}.
\end{gather}

\begin{figure}[h!]
	\begin{subfigure}{.47\textwidth} 
		{\includegraphics[width=\textwidth,height=5cm]{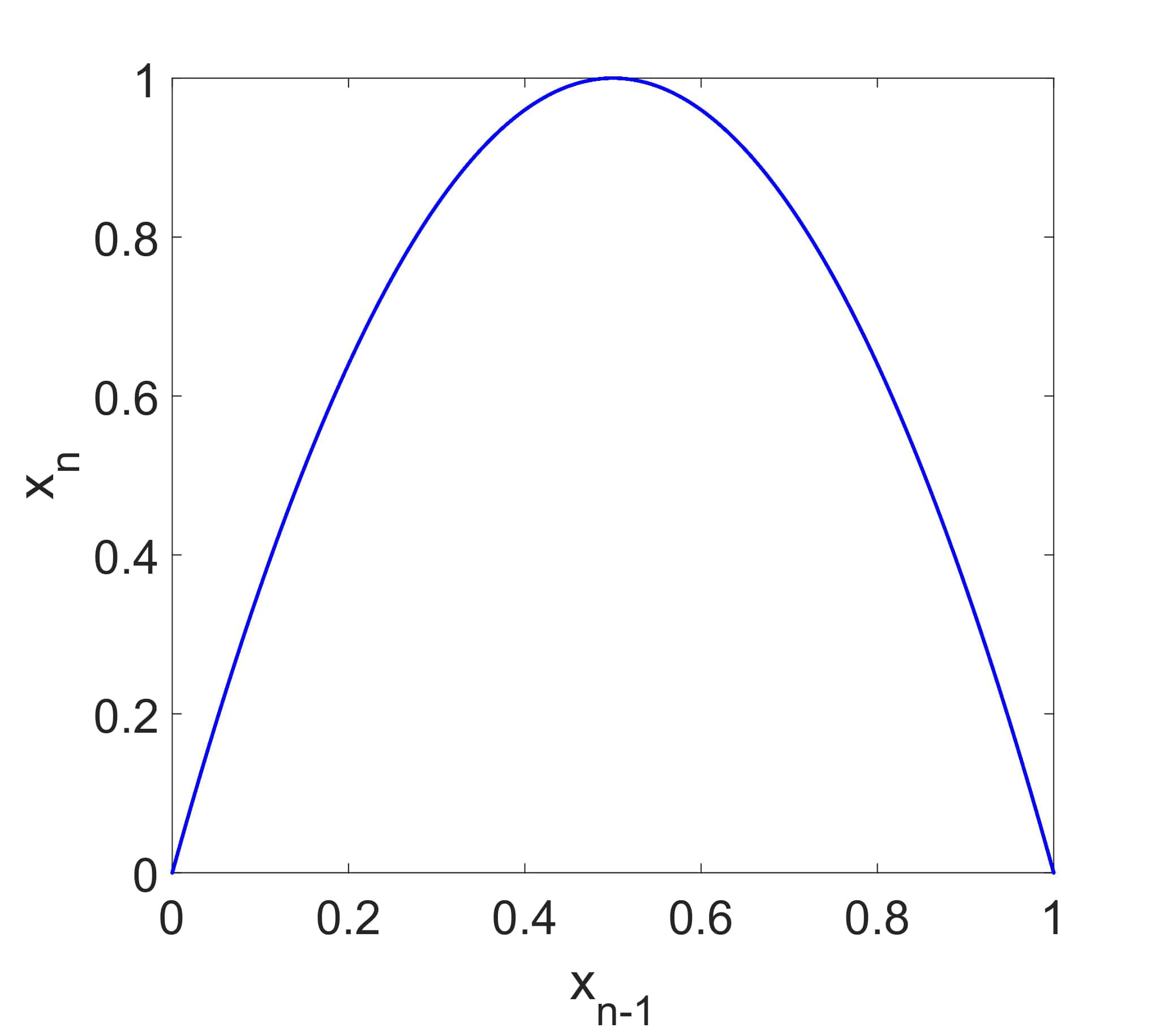}} 
		\caption{}
		\label{figLogxn}
	\end{subfigure} \hfill
	\begin{subfigure}{.47\textwidth} 
		\includegraphics[width=\textwidth,height=5cm]{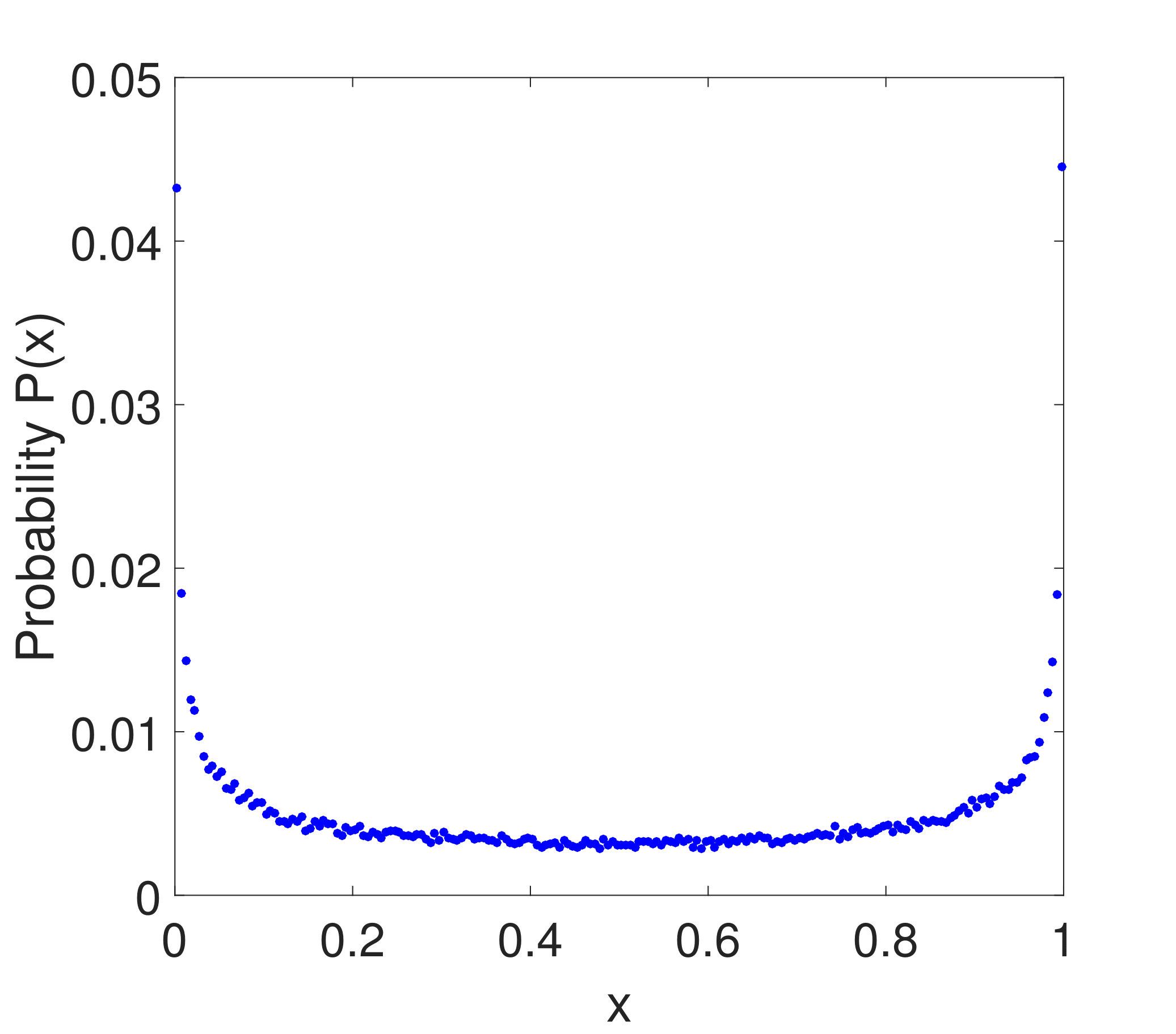} 
		\caption{}
		\label{fighistLog}
	\end{subfigure} 
	
	\begin{subfigure}{.47\textwidth} 
		{\includegraphics[width=\textwidth,height=5cm]{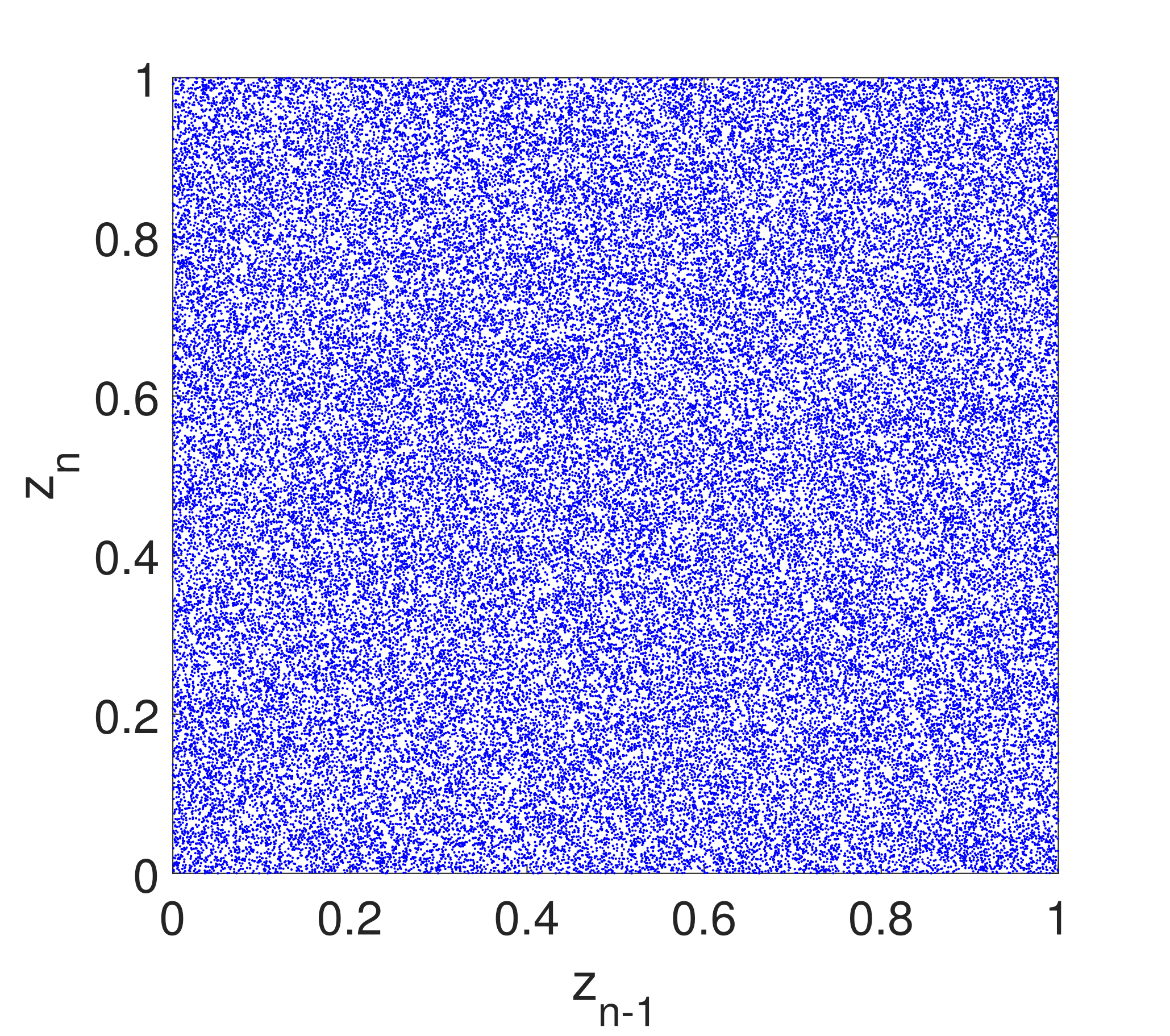}} 
		\caption{}
		\label{figzn}
	\end{subfigure}
	\hfill
	\begin{subfigure}{.47\textwidth} 
		\includegraphics[width=\textwidth,height=5cm]{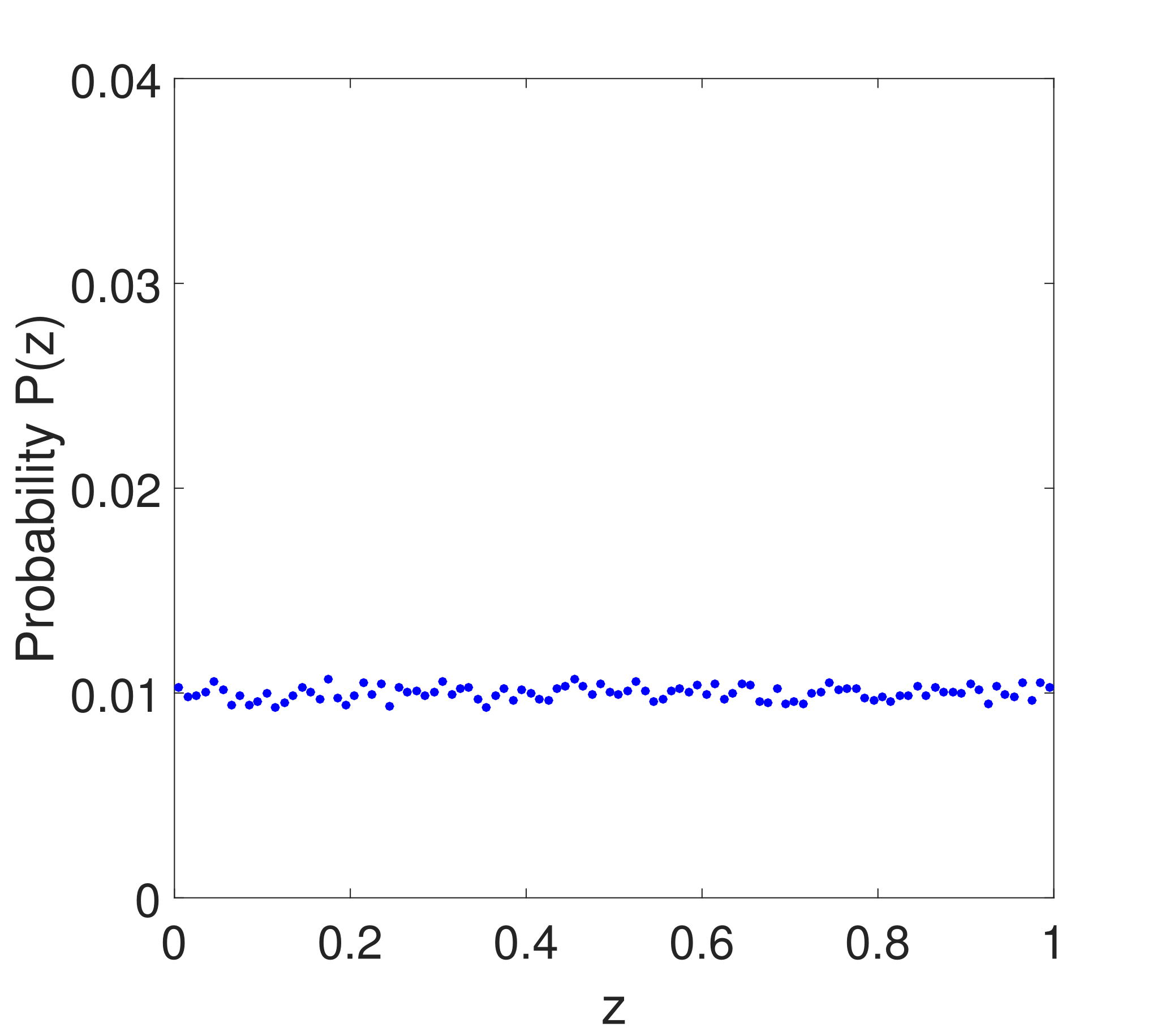}
		\caption{}
		\label{fighistZ}
	\end{subfigure}
	\caption{(a) Logistic map given by $x_n$ against $x_{n-1}$; (b) ``U-shaped'' probability distribution of the logistic map; (c) Delayed map given by $z_n$ against $z_{n-1}$; (d) Uniform probability distribution of the delayed map.} 
	\label{figLogXnZn}
\end{figure}

Finally, these time series $M1={m_{01},m_{11},m_{21},\ldots }$  and  $M2 = {m_{02}, m_{12}, m_{22}, \ldots }$ given by  \eqref{eqn:M1} and \eqref{eqn:M2}, respectively, are mixed and the operation $\mbox{mod 1}$ is applied again, this process generates a new time series $Z_i$ given as follows: 
\begin{equation}
\label{eqn:zeta}
Z_{i}=m_{i1}+m_{i2}, \mbox{mod 1}.
\end{equation}

From now on, Eq. \eqref{eqn:zeta} is referred as the delayed map.
Note that $Z_i\in[0,1)\subset \Re$.
The aim to use this approach is that with the combination of two time series with delays represented by $Z_i$, it is possible to dismiss the structure of the chaotic map used. For instance, the time series $x_n$ can reveal the map whether $x_n$ against $x_{n-1}$ is plotted as is shown in  Figure~\ref{figLogXnZn}~(a), the logistic map appears. In contrast, the time series $z_n$ can not reveal the map whether $z_n$ against $z_{n-1}$ is plotted as is shown in  Figure~\ref{figLogXnZn}~(c), the delays used are not reveled neither does the logistic map appear.  As well as, this allows us to change the characteristic ``U-shaped'' probability distribution \cite{Urias2006} by a uniform probability distribution  in the obtained time series $x_n$ and $z_n$, see Figure~\ref{figLogXnZn}~(b) and (d), respectively. This is an important characteristic that, in comparison with the chaos based schemes approach, makes easier the construction of  S-box  since all values has the same probability of occurrence in contrast with a single chaotic based schemes. 

To obtain a binary time series $s$ useful for cryptosystems, it is constructed the symbolic dynamics of $Z_i$. So the elements of $s$ are binary numbers, {\it i.e.}, $s_i (Z_i)\in \{0,1\}$. One necessary requirement for the symbolic dynamic is to obtain zeros or ones with the same probability, thus the process for getting the binary series is as follows:
\begin{equation}\label{generadorS}
s_{i}=\left\{ \begin{array}{l}

0,~~ \mbox{for} ~~0< Z_i\leq 0.5,\\
1, ~~\mbox{for} ~~0.5< Z_i< 1. \end{array}\right.
\end{equation}

A CSPRNG based on a discrete dynamical system is given from equations \eqref{logistic map} to \eqref{generadorS}.  

\subsection{The algorithm for s-box design via CSPRNG} 
\label{proposed algorithm}

In this subsection, we propose a novel design algorithm  for the creation of $n\times n$ S-boxes based on CSPRNG.  The steps of the algorithm are simple as shown below.

\begin{description}
	\item[\bf Step 1] Select initial conditions $x_{01}$ and $x_{02}$ for CSPRNG in order to generate the stream of bits $s_0, s_1, s_2,\ldots$
	\item[\bf Step 2] Generate the block sequence of $n$-bits each, $C_0=(s_0,s_{1},\ldots,s_{n-1})$, $C_1=(s_n,s_{n+1},\ldots,s_{2n-1})$, $C_2=(s_{2n},s_{2n+1},\ldots,s_{3n-1})$, $\ldots$
	\item[\bf Step 3] Convert the blocks $C_0, C_1, C_2,\ldots$ of $n$-bits to integer numbers $D_0$, $D_1$, $D_2$,$\ldots$. 
	\item[\bf Step 4]  Discard the repeated elements $D$'s to select $2^n$ different values. The rule to discard an element is as follows: if $D_i=D_j$ with $i<j$ then discard $D_j$. 
	\item[\bf Step 5] Create the S-Box with the $2^n$ different elements of $D$'s.
\end{description}

Once the procedure is over, the proposed algorithm returns a $n\times n$ S-box with distinct $2^n$ values. Note that $D_0$ is the first element of the S-box, but the second element could be not $D_1$ if $D_0=D_1$. However, has been generated the enough $2^n$ elements to build the S-box. Each block $C$'s is comprised by $n$ bits, $s_j,s_{j+1},\ldots,s_{j+n-1}$, which are related with the functions $f_i$, with $i=1,\ldots,n$. 

For example, if $n=8$, $x_{01}=0.8147$, $x_{02}=0.9058$, $\alpha_1=4$ and $\alpha_2=-2$ then the $8\times 8$ S-Box is obtained in Table \ref{S-box resesult}.
This proposed substitution box has the properties of confusion and diffusion, which are of vital importance for the block ciphers.

\begin{table}[h]
	\centering
	\resizebox*{.95\textwidth}{!}{
		\begin{tabular}{c c c c c c c c c c c c c c c c} 
			\hline
			\hline
			64 & 46 & 150 & 174 & 220 & 26 & 233 & 224 & 148 & 170 & 143 & 247 & 225 & 212 & 90 & 124 \\
			44 & 204 & 59 & 61 & 43 & 121 & 129 & 2 & 109 & 164 & 103 & 249 & 16 & 237 & 27 & 35\\
			216 & 184 & 81 & 213 & 161 & 169 & 89 & 199 & 140 & 38 & 239 & 48 & 163 & 193 & 21 & 147\\
			222 & 217 & 70 & 196 & 195 & 192 & 234 & 41 & 47 & 15 & 14 & 42 & 98 & 190 & 186 & 36\\
			242 & 51 & 60 & 87 & 24 & 104 & 189 & 55 & 118 & 111 & 231 & 120 & 8 & 226 & 7 & 141\\
			85 & 9 & 73 & 101 & 3 & 197 & 12 & 66 & 82 & 110 & 65 & 25 & 165 & 176 & 80 & 181\\
			125 & 31 & 218 & 74 & 68 & 52 & 149 & 95 & 182 & 19 & 112 & 5 & 136 & 79 & 214 & 34\\
			158 & 50 & 188 & 137 & 28 & 191 & 155 & 84 & 105 & 126 & 92 & 179 & 162 & 152 & 200 & 0\\
			171 & 142 & 240 & 203 & 88 & 160 & 32 & 202 & 99 & 18 &  100 & 97 & 145 & 53 & 194 & 93\\
			245 & 119 & 185 & 20 & 235 & 123 & 134 & 139 & 128 & 116 & 173 & 76 & 17 & 132 & 209 & 135\\
			83 & 168 & 57 & 56 & 223 & 30 & 91 & 4 & 22 & 122 & 102 & 221 & 208 & 131 & 71 & 86\\
			39 & 114 & 252 & 10 & 172 & 201 & 177 & 77 & 94 & 246 & 54 & 175 & 183 & 108 & 156 & 45\\
			219 & 210 & 40 & 130 & 113 & 153 & 13 & 166 &  58 & 23 & 253 & 215 & 238 & 33 & 198 & 248\\
			229 & 227 & 96 & 206 & 107 & 144 & 67 & 254 & 115 & 167 & 244 & 106 & 180 & 157 & 255 & 241\\
			207 & 243 & 228 & 187 & 49 & 78 & 251 & 37 & 62 & 1 & 205 & 117 & 29 & 178 & 75 & 236\\
			11 & 250 & 146 & 6 & 151 & 69 & 138 & 133 &  72 & 232 & 211 & 127 & 159 & 63 & 154 & 230\\
			\hline
			\hline
	\end{tabular}}
	\caption{The S-box generated by proposed algorithm.}
	\label{S-box resesult}
\end{table}

In the next section, is examined the performances of the proposed algorithm for the generation of S-Boxes to confirm their immunity especially against differential and linear cryptanalysis.

\section{Performance test of S-box}
\label{Performance test of S-box}
In this section we compute six important and well-known cryptographic criteria of the $8\times 8$ S-boxes. Lastly, we present  our results which are contrasted with some results presented in different published papers using other approaches. 

\subsection{Bijectivity criterion}
\label{Bijectivity criterion}
The computed value of proposed S-box is the desired value of $2^{n-1}=128$, with $n=8$, according to the formula \eqref{bijective}. So the bijectivity criterion is satisfied and the S-box proposed is a one-to-one, surjective and balance; which is a primary cryptographic criterion.

\subsection{Nonlinearity Criterion}
\label{Nonlinearity criterion}
Nonlinearity is the major requirement of any S-box design due to ensure that an S-box is not a linear function between input vectors and output vectors. The nonlinearity symbolizes the degree of dissimilarity between the Boolean function $f$ and n-bit linear function $l$. If the  function have high minimum Hamming distance is said to have high nonlinearity, {\it i.e.}, by reducing the Walsh spectrum in \eqref{nonlinearity equation}.
A S-box contains $n$ Boolean functions and  the nonlinearity of each Boolean function must be calculated. The nonlinearities of the proposed S-box are $104$, $104$, $102$, $104$, $96$, $102$, $100$ and $102$, respectively. High nonlinearity ensures the strongest ability to resist powerful modern attacks such as linear cryptanalysis.

\subsection{Strict Avalanche Criterion (SAC)}
\label{Strict Avalanche Criterion}
The Avalanche effect is used to indicate the randomness of an S-box when an input has a change. The matrix of the generated S-box can be found in Table \ref{SAC result}. For the S-box proposed, it is obtained a maximum SAC equal to $0.5781$, the minimum is $0.3906$, and its average value $0.5012$ is close to the desired value $0.5$. Based on these results, it can be concluded that the S-box generated by our proposed method fulfills the property of SAC.

\begin{table}[h]
	\centering
	\resizebox*{.7\textwidth}{!}{
		\begin{tabular}{c c c c c c c c} 
			\hline
			\hline
			
			0.5781 & 0.4844 & 0.5000 & 0.4219 & 0.4844 & 0.5156 & 0.4063 & 0.5469\\
			
			0.5156 & 0.5000 & 0.4688 & 0.5156 & 0.5469 & 0.3906 & 0.5469 & 0.4375\\
			
			0.5469 & 0.5000 & 0.5000 & 0.5469 & 0.4063 & 0.5156 & 0.4531 & 0.5313\\
			
			0.4531 & 0.5156 & 0.5000 & 0.4531 & 0.5313 & 0.5313 & 0.4844 & 0.4688\\
			
			0.5156 & 0.5469 & 0.4844 & 0.5313 & 0.5313 & 0.5625 & 0.5625 & 0.5469\\
			
			0.4063 & 0.4844 & 0.5000 & 0.4063 & 0.5625 & 0.5625 & 0.4844 & 0.5313\\
			
			0.4219 & 0.4063 & 0.5313 & 0.5313 & 0.4219 & 0.5625 & 0.4844 & 0.4844\\
			
			0.5469 & 0.5156 & 0.5469 & 0.5625 & 0.4531 & 0.5625 & 0.5781 & 0.4531\\
			
			\hline
			\hline
	\end{tabular}}
	\caption{SAC criterion result of the generated S-box.}
	\label{SAC result}
\end{table}

\subsection{Output Bits Independence Criterion (BIC)}
The BIC criterion guarantees that there is no statistic pattern or dependency between  output vectors. The BIC of the S-box generated by the proposed method is tested as described in Subsection \ref{BIC}, the results obtained are shown in Tables \ref{BIC-nonlinearity result}, \ref{BIC-SAC result} and \ref{BIC-DD result}. The mean value of BIC-nonlinearity is $103.8571$, the mean value of BIC-SAC is $0.5066$ and maximum value of DD is $8$ which indicates that S-box approximately satisfies the BIC criterion.

\begin{table}[h!]
	\centering
	\resizebox*{.5\textwidth}{!}{
		\begin{tabular}{c c c c c c c c} 
			\hline
			\hline
			0 & 104 & 104 & 106 & 104 & 106 & 106 & 102\\
			104 &   0 & 106 &  98 & 102 & 104 & 102 & 104\\
			104 & 106 &   0 & 104 & 102 &  96 & 104 & 104\\
			106 &  98 & 104 &   0 & 106 & 100 & 106 & 104\\
			104 & 102 & 102 & 106 &   0 & 102 & 100 & 102\\
			106 & 104 &  96 & 100 & 102 &   0 & 104 & 108\\
			106 & 102 & 104 & 106 & 100 & 104 &   0 & 106\\
			102 & 104 & 104 & 104 & 102 & 108 & 106 &   0\\
			\hline
			\hline
	\end{tabular}}
	\caption{BIC-nonlinearity criterion result of the generated S-box.}
	\label{BIC-nonlinearity result}
\end{table}

\begin{table}[h]
	\centering
	\resizebox*{.7\textwidth}{!}{
		\begin{tabular}{c c c c c c c c} 
			\hline
			\hline
			0 & 0.5020 & 0.5176 & 0.5137 & 0.5293 & 0.5098 & 0.4727 & 0.5059\\
			0.5020 &      0 & 0.4980 & 0.4844 & 0.5039 & 0.5313 & 0.5156 & 0.5000\\
			0.5176 & 0.4980 &      0 & 0.5039 & 0.4941 & 0.5313 & 0.5000 & 0.5020\\
			0.5137 & 0.4844 & 0.5039 &      0 & 0.5117 & 0.4980 & 0.5020 & 0.5020\\
			0.5293 & 0.5039 & 0.4941 & 0.5117 &      0 & 0.5234 & 0.5000 & 0.5137\\
			0.5098 & 0.5313 & 0.5313 & 0.4980 & 0.5234 &      0 & 0.5039 & 0.5000\\
			0.4727 & 0.5156 & 0.5000 & 0.5020 & 0.5000 & 0.5039 &      0 & 0.5156\\
			0.5059 & 0.5000 & 0.5020 & 0.5020 & 0.5137 & 0.5000 & 0.5156 &      0\\
			\hline
			\hline
	\end{tabular}}
	\caption{BIC-SAC criterion result of the generated S-box.}
	\label{BIC-SAC result}
\end{table}

\begin{table}[h]
	\centering
	\resizebox*{.4\textwidth}{!}{
		\begin{tabular}{c c c c c c c c} 
			\hline
			\hline
			0 & 2 & 6 & 2 & 4 & 6 & 4 & 6\\
			2 & 0 & 2 & 2 & 2 & 2 & 4 & 2\\
			6 & 2 & 0 & 6 & 8 & 2 & 6 & 4\\
			2 & 2 & 6 & 0 & 4 & 2 & 8 & 4\\
			4 & 2 & 8 & 4 & 0 & 2 & 0 & 2\\
			6 & 2 & 2 & 2 & 2 & 0 & 2 & 8\\
			4 & 4 & 6 & 8 & 0 & 2 & 0 & 0\\
			6 & 2 & 4 & 4 & 2 & 8 & 0 & 0\\
			\hline
			\hline
	\end{tabular}}
	\caption{The DD of the generated S-box (BIC\textendash SAC criterion).}
	\label{BIC-DD result}
\end{table}

\subsection{Criterion of equiprobable Input/Output XOR Distribution}
\label{Criterion of equiprobable input/output XOR distribution}
The equiprobable Input/Output XOR Distribution is a criterion which analyzes the effect in particular differences in input pairs of the resultant output pairs to discover the key bits. The idea is to find the high probability difference pairs for an S-Box under attack. 
The equiprobable input/output XOR distribution of generated S-box calculated by (12) is presented in Table \ref{Equiprobable input/output resesult}. Maximal value of S-box generated by the proposed method is $5$, which indicates that our S-box satisfies bound for the equiprobable Input/Output XOR Distribution criterion.

\begin{table}[h!]
	\centering
	\resizebox*{.6\textwidth}{!}{
		\begin{tabular}{c c c c c c c c c c c c c c c c}
			\hline
			\hline
			4 & 3 & 3 & 4 & 3 & 3 & 3 & 3 & 4 & 3 & 3 & 3 & 3 & 3 & 4 & 4\\
			3 & 3 & 4 & 3 & 3 & 4 & 3 & 3 & 4 & 3 & 3 & 4 & 4 & 4 & 4 & 3\\
			4 & 3 & 4 & 3 & 4 & 4 & 3 & 3 & 3 & 3 & 3 & 3 & 3 & 4 & 3 & 3\\
			4 & 3 & 3 & 3 & 4 & 4 & 4 & 4 & 3 & 4 & 5 & 4 & 3 & 2 & 3 & 3\\
			5 & 4 & 4 & 3 & 3 & 3 & 4 & 4 & 4 & 3 & 5 & 3 & 3 & 3 & 3 & 3\\
			3 & 3 & 3 & 4 & 4 & 3 & 5 & 4 & 3 & 3 & 3 & 5 & 5 & 3 & 3 & 3\\
			3 & 3 & 3 & 3 & 3 & 4 & 4 & 3 & 3 & 3 & 4 & 3 & 3 & 2 & 3 & 3\\
			3 & 2 & 3 & 3 & 3 & 4 & 3 & 3 & 3 & 3 & 3 & 4 & 3 & 3 & 3 & 3\\
			3 & 3 & 3 & 5 & 5 & 3 & 3 & 4 & 3 & 4 & 3 & 2 & 5 & 3 & 3 & 3\\
			3 & 3 & 3 & 4 & 3 & 4 & 3 & 3 & 3 & 4 & 3 & 3 & 4 & 3 & 4 & 3\\
			4 & 3 & 4 & 3 & 2 & 3 & 3 & 4 & 3 & 3 & 3 & 3 & 3 & 4 & 3 & 3\\
			3 & 4 & 3 & 3 & 3 & 3 & 3 & 3 & 3 & 4 & 3 & 3 & 3 & 3 & 4 & 4\\
			3 & 3 & 3 & 3 & 3 & 4 & 3 & 3 & 2 & 4 & 3 & 3 & 4 & 4 & 3 & 3\\
			4 & 3 & 4 & 3 & 4 & 4 & 3 & 4 & 4 & 3 & 4 & 4 & 3 & 3 & 3 & 3\\
			3 & 4 & 3 & 3 & 3 & 3 & 3 & 3 & 3 & 4 & 4 & 3 & 3 & 3 & 3 & 3\\
			3 & 3 & 3 & 5 & 4 & 5 & 4 & 3 & 3 & 5 & 3 & 3 & 4 & 3 & 5 & -\\
			\hline
			\hline
	\end{tabular}}
	\caption{Equiprobable Input/Output XOR Distribution approach table for the generated S-box.}
	\label{Equiprobable input/output resesult}
\end{table}

\subsection{MELP criterion of the generated S-box}
\label{Melp Criterion}
The MELP value is calculated via linear approximations to model nonlinear steps. The final goal is to recover the key bits or part of the key bits. MELP studies the statistical correlation between the input and the output.
This criterion of proposed S-Box is computed according to equation \eqref{MEDP equation} and the average  value is $0.0176$.

\subsection{Comparative results}
\label{Comparative results}

In this section a performance comparison of the S-box created using the algorithm described in Subsection \ref{proposed algorithm} is presented. In the Table \ref{Comparison} the criteria values for our S-box and a set of widely known boxes (standard and chaos based S-box) are shown. From this Table \ref{Comparison}, it can be seen that the generated S-box fulfills the most important condition, bijectivity, and accomplish a good similarity to the rest of the test values expected \cite{jakimoski2001chaos,belazi2017efficient, ccavucsouglu2017novel,cui2007new,tran2008gray,daemen1991design,tang2005method,khan2012novel, belazi2017chaos,ul2017designing,ozkaynakconstruction}. Mainly it shows better performance in the tests related to attacks (MELP and equiprobable Input/Output XOR Distribution). Also is a methodology based on a system with simple operations that generates sequences with complex behavior.

\begin{table}[h!]
	\centering
	\resizebox*{1.00\textwidth}{!}{
		\begin{tabular}{|c|c|ccc|ccc|ccc|c|c|}
			\hline
			\multirow{2}{*}{} & \multicolumn{1}{l|}{\multirow{2}{*}{\textbf{Bijective}}} & \multicolumn{3}{c|}{\textbf{Nonlinearity}} & \multicolumn{3}{c|}{\textbf{SAC}} &  \multicolumn{3}{c|}{\textbf{BIC}} & \multirow{2}{*}{\textbf{I/O XOR}} & \multirow{2}{*}{\textbf{MELP}} \\
			& \multicolumn{1}{l|}{} & \textbf{min} & \textbf{max} & \textbf{avg} & \textbf{min} & \textbf{max} & \textbf{avg} & \textbf{SAC} &\textbf{Nonlinearity} &\textbf{DD} &    &  \\ \hline
			Skipjack S-Box  \cite{hussain2011analyses}&  123  &  100 & 108 & 105.1250&    0.3906&    0.5938 &   0.5027 &      0.5003 & 104.03 & 109 &  0.0469  &  0.0137\\
			APA S-Box \cite{cui2007new} & 128 & 112 & 112  &112  & 0.4375 &  0.5625 &   0.5007 & 0.4997 & 112 &      112 & 0.0156 &   0.0039\\
			Gray S-Box \cite{tran2008gray} & 128 & 112 &  112 &  112 & 0.4375 &   0.5625 &   0.4998  & 0.5026 &  112      & 112 &  0.0156  &  0.0039 \\
			AES S-Box \cite{daemen1991design} & 128 &      112 & 112 & 112 &  0.4531 & 0.5625  &  0.5049  &  0.5046 & 112  &     112 &     0.0156  &  0.0039 \\
			Ref. \cite{jakimoski2001chaos} & 128 & 98 & 107 & 103.25 & 0.3828 & 0.5938 & 0.5059 & 0.5033 & 104.21 & 108 & 0.0469 & 0.0166 \\
			Ref. \cite{tang2005method} & 129 & 103 & 109 & 104.875 & 0.3984 & 0.5703 & 0.4966 & 0.5044 & 102.96& 109 & 0.0391 & 0.0176 \\
			Ref. \cite{khan2012novel} & 128 & 96 & 106 & 103 & 0.3906 & 0.6250 & 0.5039 & 0.5010 & 100.35 & 106 & 0.5000 & 0.0220 \\
			Ref. \cite{belazi2017chaos} & 128 & 112 & 112 & 112 & 0.4219 & 0.5469 & 0.5115 & 0.4982 & 108.71& 112 & 0.0313 & 0.0120 \\
			Ref. \cite{belazi2017efficient} & 128 & 102 & 108 & 105.25 & 0.4375 & 0.5781 & 0.5056 & 0.5019 & 103.78 & 108 & 0.0391 & 0.0244 \\
			Ref. \cite{ccavucsouglu2017novel} & 128 & 104 & 110 & 106.25 & 0.4219 & 0.5938 & 0.5039 & 0.5059 & 103.35 & 108 & 0.0391 & 0.0198 \\
			Ref. \cite{ul2017designing} & 128 & 102 & 108 & 106 & 0.4219 & 0.5938 & 0.5002 & 0.5016 & 104.42 & 108 & 0.0391 & 0.0220 \\
			Ref. \cite{ozkaynakconstruction}-1 & 128 & 106 & 108 & 106.75 & 0.3906 & 0.6094 & 0.4941 & 0.5013 & 104.28 & 108 & 0.0391 & 0.0156 \\
			Ref. \cite{ozkaynakconstruction}-2 & 128 & 106 & 108 & 106.75 & 0.4063 & 0.5938 & 0.4971 & 0.5008 & 102.92 & 106 & 0.0391 & 0.0198 \\
			The proposed S-box & 128 & 96 & 104 & 101.75 & 0.3906 & 0.5781 & 0.5012 & 0.5066 & 103.42 & 108 & 0.0391 & 0.0176 \\ \hline
		\end{tabular}
	}
	\caption{Comparison of recent chaos-base designed S-Boxes and S-boxes used in typical block ciphers.}
	\label{Comparison}
\end{table}

\section{Dynamical generation of S-boxes and its application}
\label{Application of dynamical S-boxes}

\begin{figure}[h!]
	\centering
	\begin{subfigure}{.37\textwidth} 
		\includegraphics[width=\textwidth,height=5cm]{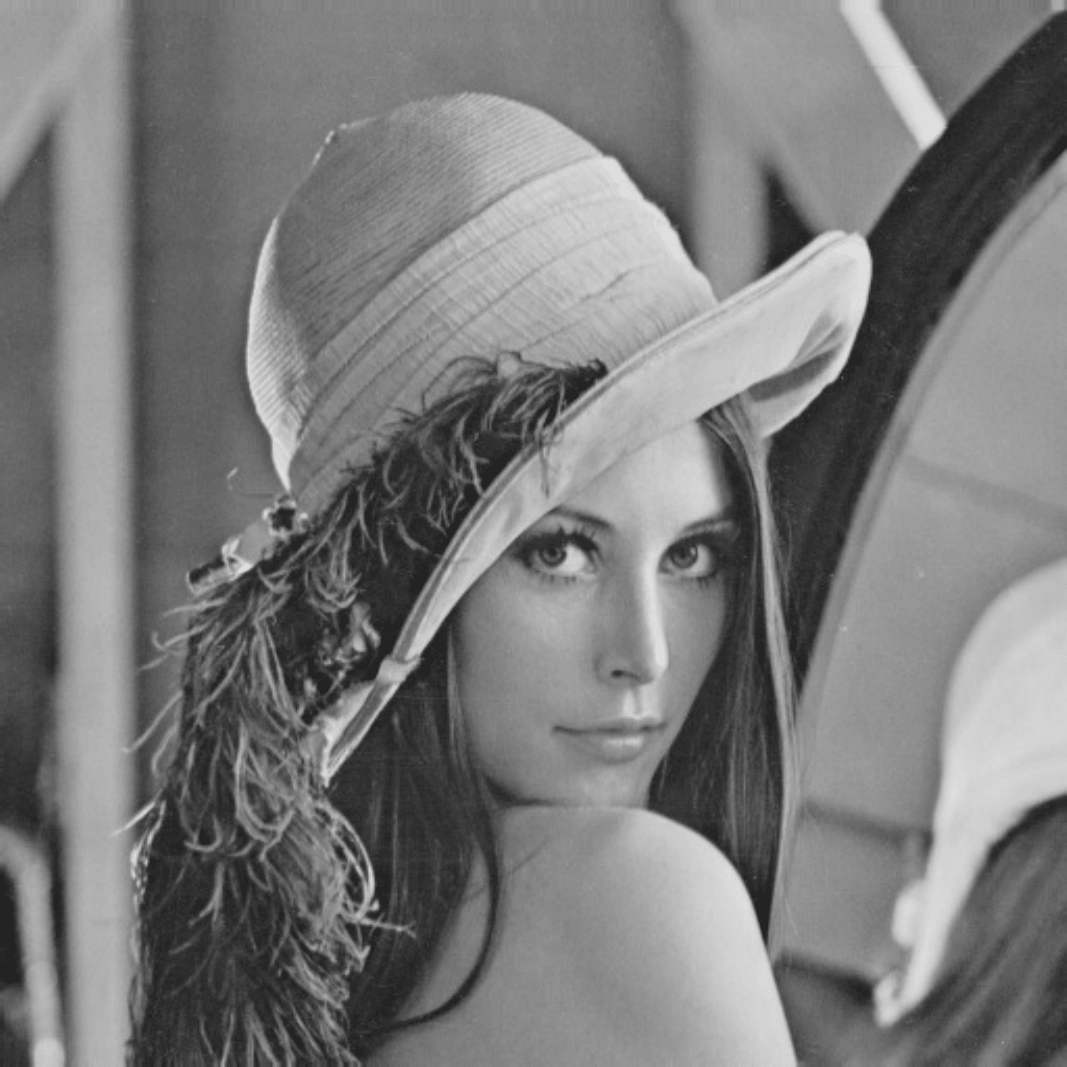}
		\caption{}
		\label{fig:sfig1}
	\end{subfigure} ~
	\begin{subfigure}{.4\textwidth} 
		\includegraphics[width=\textwidth,height=5cm]{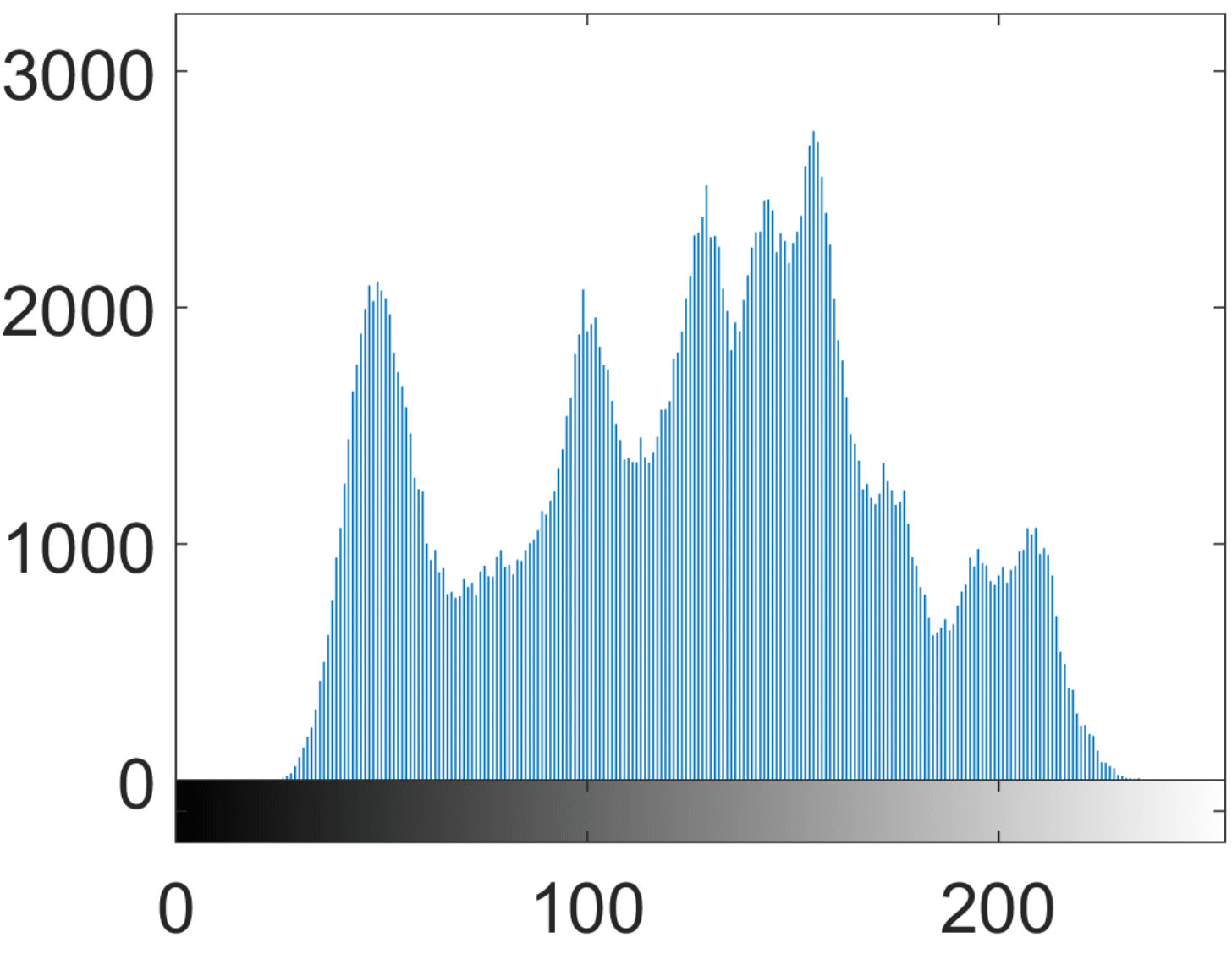}
		\caption{}
		\label{fig:sfig2}
	\end{subfigure}
	
	\begin{subfigure}{.37\textwidth} 
		{\includegraphics[width=\textwidth,height=5cm]{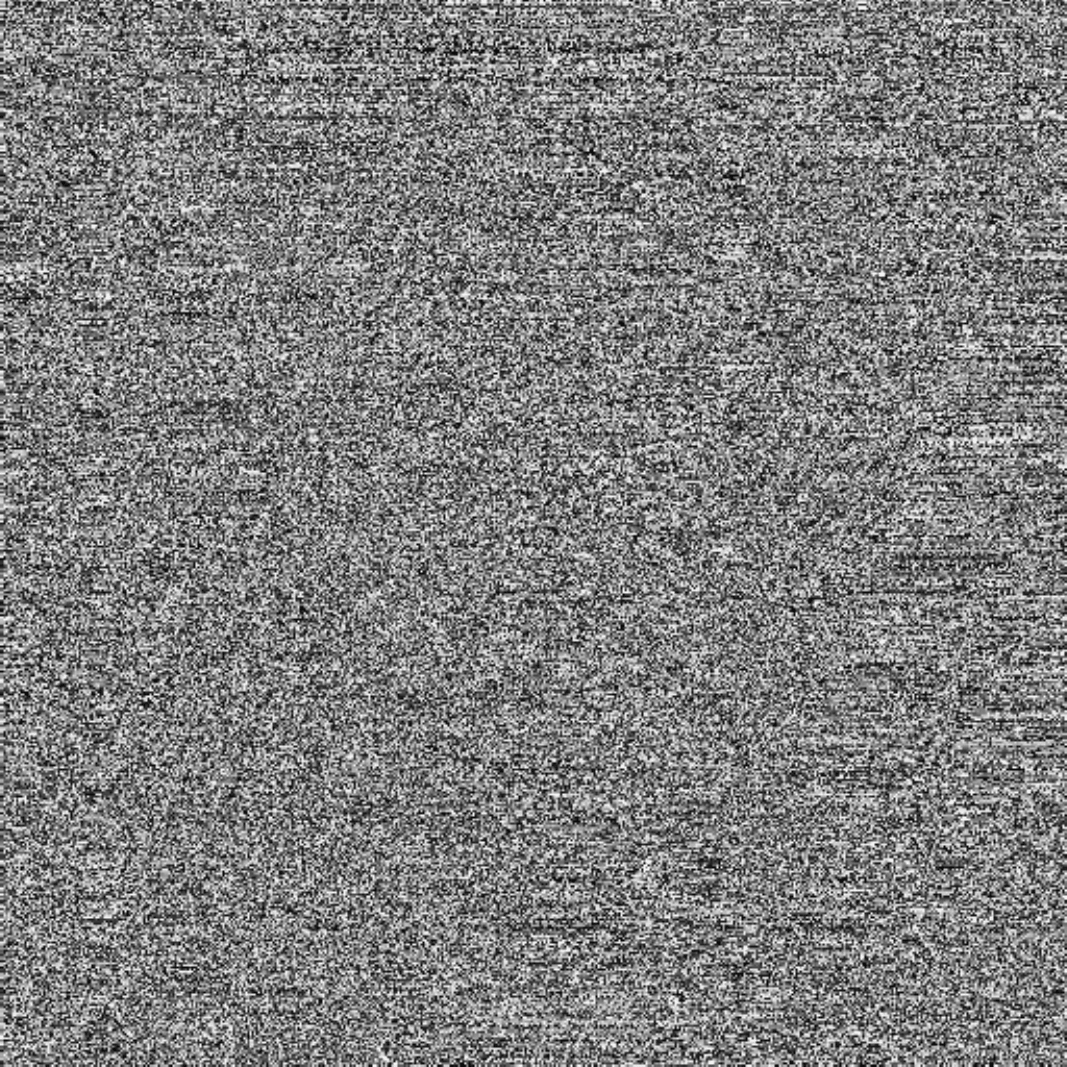}}
		\caption{}
		\label{fig:sfig3}
	\end{subfigure} ~
	\begin{subfigure}{.4\textwidth} 
		{\includegraphics[width=\textwidth,height=5cm]{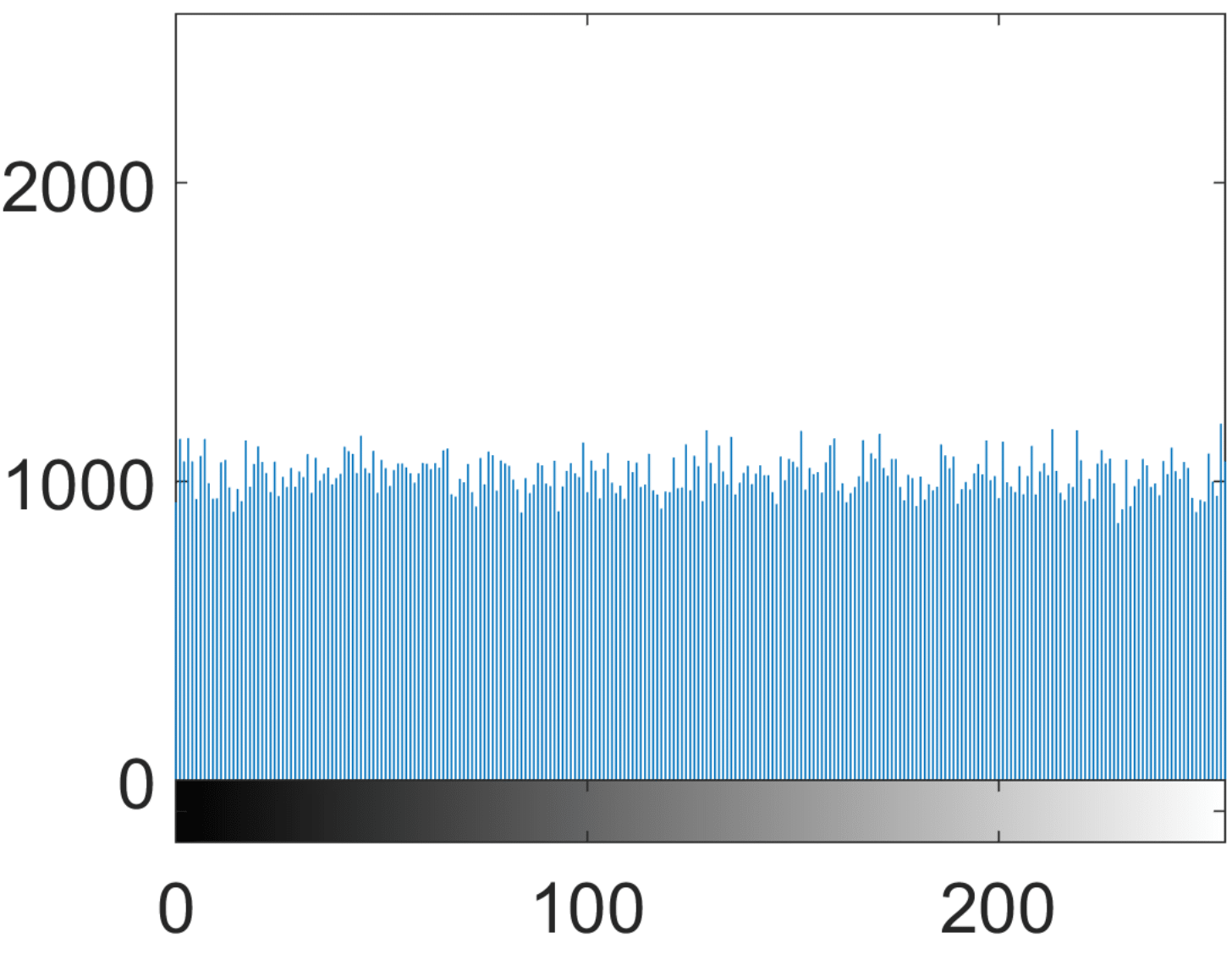}}
		\caption{}
		\label{fig:sfig4}
	\end{subfigure}
	
	\caption{Plain image of Lena, encrypted image and their histograms.} 
	\label{fig:Lena}
\end{figure}

The Alberti cipher was one of the first polyalphabetic ciphers where the principle is based on substitution, using multiple substitution alphabets such that the output has a uniform distribution.
Taking this idea of polyalphabetic ciphers, it is presented an application of dynamical generation of S-boxes, {\it i.e.}, a particular intensity of a pixel given can be substituted by different intensities in the same round. 
Usually an S-box is used to substitute all the pixels of an image of size $p\times q$ in the same way. The idea of polyalphabetic ciphers is to use a dynamical S-box to achieve this aim looking for a uniform distribution. For instance, our dynamical S-box belongs to a class of S-boxes given by $p$ elements (S-boxes) generated by the algorithm presented in the Subsection~\ref{proposed algorithm}.
The approach to get uniform distribution is given by applying dynamical S-box which changes with each pixel row, {\it i.e.}, the dynamical S-box is modified with a different S-box. Figure \ref{fig:Lena} shows the original Lena image and the codified Lena image alone with their gray scale pixels distribution. 

In crypthography a uniform distribution is always desired, since this property  was achieved by simple substitution of our S-boxes, a good result is spected for a full cryptographic algorithm based on our S-boxes.

It is important to point out that this is not an encryption algorithm, but a simple and useful approach intended to bring to light possible applications of our dynamical S-boxes.

\section{\label{sec:Conclusions}Concluding remarks}
\label{conclusions}

In this work, a simple algorithm was proposed to generate $n\times n$ S-boxes by means of using a pseudo-random bit generator (PRBG) based on two lag time series of the logistic map. The mixed of this two time series favors a uniform distribution in addition to hide the chaotic map used. To evaluate performance of the proposed S-box, several statistical tests were carried out. The results of numerical analysis of these cryptographic strong S-box generated by the algorithm have also shown that all criteria for a good S-box were fulfilled and with high immunity to resist differential cryptanalysis and linear cryptanalysis. The performance test result was compared with other S-Boxes which were reported in the literature. Finally, an application based on simple and useful approach to bring uniform distribution was presented.

\section*{Acknowledgements}
B.B.C.Q  is a doctoral fellow of the CONACYT in the Postgraduate Program in control and dynamical systems of the DMAp-IPICYT. 




\begin{thebibliography}{10}
	\expandafter\ifx\csname url\endcsname\relax
	\def\url#1{\texttt{#1}}\fi
	\expandafter\ifx\csname urlprefix\endcsname\relax\def\urlprefix{URL }\fi
	\expandafter\ifx\csname href\endcsname\relax
	\def\href#1#2{#2} \def\path#1{#1}\fi
	
	\bibitem{Shannon1949}
	C.~E. Shannon, Communication theory of secrecy systems, Bell System Technical
	Journal 28~(4) (1949) 656--715.
	
	\bibitem{Adams1990structured}
	C.~Adams, S.~Tavares, The structured design of cryptographically good s-boxes,
	Journal of Cryptology 3~(1) (1990) 27--41.
	
	\bibitem{jakimoski2001chaos}
	G.~Jakimoski, L.~Kocarev, Chaos and cryptography: block encryption ciphers
	based on chaotic maps, IEEE Transactions on Circuits and Systems I:
	Fundamental Theory and Applications 48~(2) (2001) 163--169.
	
	\bibitem{CHEN}
	G.~Chen, A novel heuristic method for obtaining s-boxes, Chaos, Solitons \&
	Fractals 36~(4) (2008) 1028 -- 1036.
	
	\bibitem{WANG}
	Y.~Wang, K.-W. Wong, X.~Liao, T.~Xiang, A block cipher with dynamic s-boxes
	based on tent map, Communications in Nonlinear Science and Numerical
	Simulation 14~(7) (2009) 3089 -- 3099.
	
	\bibitem{LAMBIC}
	D.~Lambi\'{c}, A novel method of s-box design based on chaotic map and
	composition method, Chaos, Solitons \& Fractals 58 (2014) 16 -- 21.
	
	\bibitem{belazi2017efficient}
	A.~Belazi, M.~Khan, A.~A.~A. El-Latif, S.~Belghith, Efficient cryptosystem
	approaches: S-boxes and permutation--substitution-based encryption, Nonlinear
	Dynamics 87~(1) (2017) 337--361.
	
	\bibitem{OZKAYNAK}
	F.~{\"O}zkaynak, A.~B. {\"O}zer, A method for designing strong s-boxes based on
	chaotic {Lorenz} system, Physics Letters A 374~(36) (2010) 3733 -- 3738.
	
	\bibitem{Liu}
	G.~Liu, W.~Yang, W.~Liu, Y.~Dai, Designing s-boxes based on 3-d four-wing
	autonomous chaotic system, Nonlinear Dynamics 82~(4) (2015) 1867--1877.
	
	\bibitem{ccavucsouglu2017novel}
	{\"U}.~{\c{C}}avu{\c{s}}o{\u{g}}lu, A.~Zengin, I.~Pehlivan, S.~Ka{\c{c}}ar, A
	novel approach for strong s-box generation algorithm design based on chaotic
	scaled {Zhongtang} system, Nonlinear Dynamics 87~(2) (2017) 1081--1094.
	
	\bibitem{Guesmi}
	R.~Guesmi, M.~A.~B. Farah, A.~Kachouri, M.~Samet, A novel design of chaos based
	s-boxes using genetic algorithm techniques, in: 2014 IEEE/ACS 11th
	International Conference on Computer Systems and Applications (AICCSA), 2014,
	pp. 678--684.
	
	\bibitem{Tian}
	Y.~Tian, Z.~Lu, S-box: Six-dimensional compound hyperchaotic map and artificial
	bee colony algorithm, Journal of Systems Engineering and Electronics 27~(1)
	(2016) 232--241.
	
	\bibitem{ALVAREZ}
	G.~Alvarez, S.~Li, Some basic cryptographic requirements for chaos-based
	cryptosystems, International Journal of Bifurcation and Chaos 16~(08) (2006)
	2129--2151.
	
	\bibitem{Azkaynak2013}
	F.~{\"O}zkaynak, S.~Yavuz, Designing chaotic s-boxes based on time-delay
	chaotic system, Nonlinear Dynamics 74~(3) (2013) 551--557.
	
	\bibitem{ZHOU}
	Y.~Zhou, L.~Bao, C.~P. Chen, Image encryption using a new parametric switching
	chaotic system, Signal Processing 93~(11) (2013) 3039 -- 3052.
	
	\bibitem{Garcia2014Campos}
	M.~Garc\'ia-Mart\'inez, E.~Campos-Cant\'on, Pseudo-random bit
	generator based on lag time series, International Journal of Modern Physics C
	25~(04) (2014) 1350105.
	
	\bibitem{Dimitris}
	D.~Souravliasa, K.~E. Parsopoulos, G.~C. Meletiou, Designing bijective s-boxes
	using algorithm portfolios with limitedtime budgets, Applied Soft Computing
	59~(1) (2017) 475--486.
	
	\bibitem{adams1989good}
	C.~Adams, S.~Tavares, Good s-boxes are easy to find, in: G.~Brassard (Ed.),
	Advances in Cryptology --- CRYPTO' 89 Proceedings, Springer New York, New
	York, NY, 1990, pp. 612--615.
	
	\bibitem{Tian2017Chaos}
	Y.~Tian, Z.~Lu, Chaotic s-box: Intertwining logistic map and bacterial foraging
	optimization, Mathematical Problems in Engineering (2017) 1--11.
	
	\bibitem{Millan1998}
	W.~Millan, How to improve the nonlinearity of bijective s-boxes, in: C.~Boyd,
	E.~Dawson (Eds.), Information Security and Privacy, Springer Berlin
	Heidelberg, Berlin, Heidelberg, 1998, pp. 181--192.
	
	\bibitem{Webster1986}
	A.~F. Webster, S.~E. Tavares, On the design of s-boxes, in: H.~C. Williams
	(Ed.), Advances in Cryptology --- CRYPTO '85 Proceedings, Springer Berlin
	Heidelberg, Berlin, Heidelberg, 1986, pp. 523--534.
	
	\bibitem{biham1991differential}
	E.~Biham, A.~Shamir, Differential cryptanalysis of {DES}-like cryptosystems,
	in: A.~J. Menezes, S.~A. Vanstone (Eds.), Advances in Cryptology-CRYPT0' 90,
	Springer Berlin Heidelberg, Berlin, Heidelberg, 1991, pp. 2--21.
	
	\bibitem{may1976simple}
	R.~M. May, Simple mathematical models with very complicated dynamics, Nature
	261~(5560) (1976) 459.
	
	\bibitem{dendrinos1993socio}
	D.~S. Dendrinos, M.~Sonis, Socio-spatial stocks and antistocks; the logistic
	map in real space, The Annals of Regional Science 27~(4) (1993) 297--313.
	
	\bibitem{devaneyintroduction}
	R.~L. Devaney, An introduction to chaotic dynamical systems, Westview Press.
	
	\bibitem{Lyapunovexponents2004}
	C.~Li, G.~Chen, Estimating the lyapunov exponents of discrete systems, Chaos
	14~(2) (2004) 343--346.
	
	\bibitem{Yang2012}
	C.~Yang, C.~Q. Wu, P.~Zhang, Estimation of lyapunov exponents from a time
	series for n-dimensional state space using nonlinear mapping, Nonlinear
	Dynamics 69~(4) (2012) 1493--1507.
	
	\bibitem{Urias2006}
	J.~Ur{\'i}as, E.~Campos, N.~F. Rulkov, Random Finite Approximations of Chaotic
	Maps, Springer, New York, NY, 2006, pp. 231--242.
	
	\bibitem{cui2007new}
	L.~Cui, Y.~Cao, A new s-box structure named affine-power-affine, International
	Journal of Innovative Computing, Information and Control 3~(3) (2007)
	751--759.
	
	\bibitem{tran2008gray}
	M.~T. Tran, D.~K. Bui, A.~D. Duong, Gray s-box for advanced encryption
	standard, in: 2008 International Conference on Computational Intelligence and
	Security, Vol.~1, IEEE, 2008, pp. 253--258.
	
	\bibitem{daemen1991design}
	J.~Daemen, V.~Rijmen, {AES} proposal: Rijndael, available:
	http://csrc.nist.gov/archive/aes/rijndael/Rijndael-ammended.pdf (1999).
	
	\bibitem{tang2005method}
	G.~Tang, X.~Liao, A method for designing dynamical s-boxes based on discretized
	chaotic map, Chaos, Solitons \& Fractals 23~(5) (2005) 1901--1909.
	
	\bibitem{khan2012novel}
	M.~Khan, T.~Shah, H.~Mahmood, M.~A. Gondal, I.~Hussain, A novel technique for
	the construction of strong s-boxes based on chaotic {Lorenz} systems,
	Nonlinear Dynamics 70~(3) (2012) 2303--2311.
	
	\bibitem{belazi2017chaos}
	A.~Belazi, A.~A.~A. El-Latif, A.-V. Diaconu, R.~Rhouma, S.~Belghith,
	Chaos-based partial image encryption scheme based on linear fractional and
	lifting wavelet transforms, Optics and Lasers in Engineering 88 (2017)
	37--50.
	
	\bibitem{ul2017designing}
	F.~ul~Islam, G.~Liu, Designing s-box based on {4D}-4wing hyperchaotic system,
	3D Research 8~(1) (2017) 1--9.
	
	\bibitem{ozkaynakconstruction}
	F.~{\"O}zkaynak, Construction of robust substitution boxes based on chaotic
	systems, Neural Computing and Applications (2017) 1--10.
	
	\bibitem{hussain2011analyses}
	I.~Hussain, T.~Shah, M.~A. Gondal, Y.~Wang, Analyses of {SKIPJACK} s-box, World
	Appl. Sci. J. 13~(11) (2011) 2385--2388.
	
\end{thebibliography}
\end{document}